\documentclass[numbers,webpdf,imanum]{ima_authoring_template}%
\usepackage{amsfonts}
\usepackage{amsmath}
\usepackage{amssymb}
\usepackage{graphicx}
\usepackage{subfig}
\usepackage{booktabs}
\usepackage{bm}
\usepackage{xcolor}
\usepackage{natbib}
\setcitestyle{authoryear,round,citesep={;\!}}
\graphicspath{{figs/}}
\catcode`\@=11
\numberwithin{equation}{section}
\numberwithin{figure}{section}
\numberwithin{table}{section}
\numberwithin{equation}{section}
\newtheorem{thm}{Theorem}[section] 

\newtheorem{lem}[thm]{Lemma}
\newtheorem{rem}[thm]{Remark}
\newtheorem{corollary}[thm]{Corollary}
\graphicspath{{figs/}}
\usepackage{bm}
\usepackage{longtable}
\usepackage{float}
\usepackage{hyperref}

\makeatletter
\newcommand{\removelatexerror}{\let\@latex@error\@gobble}
\makeatother
\usepackage{listings}%
\usepackage{makecell}

\newcommand{\mbR}{\mathbb{R}}
\newcommand{\td}{\mathrm{d}}
\newcommand{\htheta}{\hat{\theta}}

\begin{document}
\DOI{}
\copyrightyear{}
\vol{}
\pubyear{}
\access{}
\appnotes{}
\copyrightstatement{}
\firstpage{1}

\title[general transformation neural networks]{General transformation neural networks: A class of parametrized functions for high-dimensional function approximation}

\author{Xiaoyang Wang
\address{\orgdiv{School of Mathematical Sciences}, \orgname{University of Electronic Science and Technology of China}, \orgaddress{\state{Sichuan}, \postcode{610054} \country{China}}}}
\author{Yiqi Gu*
\address{\orgdiv{School of Mathematical Sciences}, \orgname{University of Electronic Science and Technology of China}, \orgaddress{\state{Sichuan}, \postcode{610054}, \country{China}}}}
\authormark{X. Wang and Y. Gu}
\corresp[*]{Corresponding author: yiqigu@uestc.edu.cn}

\received{Date}{0}{Year}
\revised{Date}{0}{Year}
\accepted{Date}{0}{Year}


\abstract{We propose a novel class of neural network-like parametrized functions, i.e., general transformation neural networks (GTNNs), for high-dimensional approximation. Conventional deep neural networks sometimes perform less accurately on learning problems trained with gradient descent, especially when the target function is oscillatory. To improve accuracy, we generalize the neuron's affine transformation to a broader class of functions that can capture complex shapes and offer greater capacity. Specifically, we discuss three types of GTNNs in detail: the cubic, quadratic and trigonometric transformation neural networks (CTNNs, QTNNs and TTNNs). We perform an approximation error analysis of GTNNs, presenting their universal approximation properties for continuous functions, error bounds for Barron-type functions and error bounds of deep architectures. Several numerical examples of regression problems are presented, demonstrating that CTNNs/QTNNs/TTNNs achieve higher accuracy than conventional fully connected neural networks.}

\allowdisplaybreaks
\sloppy
\maketitle

\section{Introduction}  
Deep neural networks (DNNs) are parametrized functions that serve as good approximators of high-dimensional functions, making them widely used in artificial intelligence, data science, and applied mathematics, among other fields. There are research studies of the approximation theory of DNNs, showing particular advantages over traditional parametrized functions. The research traces back to the pioneering work \citep{Barron1993,Cybenko1989,Hornik1989}, on the universal approximation of shallow networks with sigmoid activation functions. The recent research focus was on the approximation rate of DNNs for various function spaces in terms of the number of network parameters, e.g., Barron functions \citep{E2019,E2022}, smooth functions \citep{E2018,Liang2017,Lu2017,Montanelli2019,Suzuki2019,Yarotsky2017,Lu2021}, piecewise smooth functions \citep{Petersen2018}, band-limited functions \citep{Montanelli2021}, continuous functions \citep{Shen2020,Shen2019,Yarotsky2018}. For some types of target functions, the approximation error of DNNs grows with the dimension at most polynomially; namely, DNNs can lessen or overcome the curse of dimensionality in high-dimensional approximation \citep{Chui2018,Guliyev2018,Maiorov1999}. Comparatively, linear structures such as polynomials, Fourier series, and finite elements suffer from the curse of dimensionality in approximation \citep{Chkifa2015,Parhi2023}, whose error orders are inversely proportional to the dimension for common target functions. 

\subsection{Motivation}\label{Motivation}
However, conventional DNNs are less efficient in approximating some special functions. Recent research has shown that DNNs exhibit spectral biases: during gradient-based training, they tend to initially capture the lower frequencies of the target function while overlooking higher ones. The spectral bias and the frequency principle have been extensively observed through empirical works \citep{Xu2019,Zhi2020,Cao2021,Rahaman2019,Xu2024} and supported by theoretical analyses \citep{Zhi2020,Cao2021,Luo2019,Luo2022,Zhang2021,Bordelon2020,Luo2022_2}. Therefore, it is difficult for DNNs to learn high-frequency functions, oscillatory functions and singular functions, particularly when the dimension is high. 

Comparatively, finite element functions can handle high-frequency functions effectively using dense meshes. On each mesh element, the target function behaves smoothly, allowing it to be well approximated by low-order elements such as linear elements. To achieve higher accuracy, one can further use spectral element functions, which employ a linear combination of high-degree polynomials for approximation on each mesh element. In other words, spectral element functions provide a richer class of shape functions than finite elements, hence having a greater representational capability.

\subsection{Our contribution}
In this work, inspired by the above discussion of finite/spectral element functions, we propose a novel class of DNN-like parametrized functions. Note that DNNs with ReLU activation are exactly piecewise linear functions, which are composed of multiple layers of abstract neurons $h(x)=\sigma(Wx+b)$, where $W$ and $b$ are the weight matrix and bias vector, respectively, and $\sigma$ is the ReLU activation. These affine transformations $x\mapsto Wx+b$ applied by ReLU play the role of high-dimensional domain partitioning, and the ReLU DNNs can be analogous to high-dimensional finite elements. We now aim to generalize the linear shape function to more complex functions, such as high-degree polynomials or trigonometric functions. Specifically, we multiply the affine term $Wx+b$ by a specific vector-valued parametrized function $\psi(x;\theta)$, obtaining a general neuron $g(x)=\sigma((Wx+b)\ast\psi(x;\theta))$, where $\ast$ means the entry-wise multiplication. These functions $\psi(x;\theta)$ can be chosen as high-dimensional polynomials or other types of functions, which act as shape functions on the mesh elements. Based on this, we propose the general transformation neural networks (GTNNs), which are composite functions of multi-layer general neurons and are heuristically similar to high-dimensional spectral elements. Briefly, the GTNN is a generalization of the conventional DNN; the part acted upon by the activation function is a general parametrized function rather than a parametrized linear function.

We also analyze the approximation properties of shallow and deep GTNNs and, in particular, provide detailed discussions of three representative GTNN architectures: the quadratic transformation neural network (QTNN), the cubic transformation neural network (CTNN), and the trigonometric transformation neural network (TTNN). We first investigate the universal approximation property, which indicates that the two-layer GTNN acts as an approximator of specific functions to any desired degree of accuracy. Then, we derive error bounds that do not depend on the dimension explicitly for approximating Barron-type functions, so the approximation lessens or overcomes the curse of dimensionality in some sense. We also derive an error bound of deep architectures in terms of the depth and width of the GTNNs. Finally, numerical experiments on the least squares regression problem are conducted, demonstrating the advantages of QTNNs/CTNNs/TTNNs over conventional DNNs.  

\subsection{Related works and comparison}
One recent line of research on network architectures concerns quadratic neural networks. In previous studies \citep{Fan2018}, quadratic neurons and deep quadratic networks have been proposed. A quadratic neuron with a single Hadamard (element-wise) interaction computes the product, thereby causing the input $x$ to enter the nonlinear interaction twice (i.e., quadratic), and the two-layer quadratic networks are a linear combination of quadratic neurons. The work in \citep{Fan2020} also provably demonstrates the global universal approximation of a quadratic network with ReLU activation function. For quadratic deep networks, the general backpropagation algorithm enables the network training process \citep{Fan2017}. Evidence in \citep{Fan2025,Fan2023} indicates that the performance gains of quadratic networks are not merely a consequence of having more parameters but arise from greater intrinsic expressivity: quadratic units natively capture nonlinear pairwise interactions that are cumbersome for conventional neurons. The study in \citep{Bu2021} proposed quadratic residual networks that incorporate quadratic residual terms before activation functions, demonstrating parameter efficiency in solving both forward and inverse PDE problems. For image classification tasks, the work in \citep{Mantini2020} developed convolutional quadratic neural networks, which replace linear neurons with quadratic ones in convolutional layers, significantly improving classification accuracy with exponential linear unit activations. Further enhancing architectural efficiency, the authors in \citep{Jiang2020} introduced nonlinear convolutional networks by inserting quadratic convolutional units to map features to higher-dimensional spaces and employed a genetic algorithm-based training to mitigate computational costs.

Another type of research focuses on trainable activation functions. Generally, one can relax the rigidity of fixed nonlinear activations by learning their shape from data, thereby improving expressivity without substantially increasing depth. Classic parametric rectifiers tune a few scalars to control slope or curvature (e.g., PReLU/PELU/SReLU) and have been shown to alleviate vanishing gradients and accelerate convergence \citep{He2015,Trottier2017,Jin2016}.

Beyond previous works, we propose a more general type of neuron $g(x)=\sigma((Wx+b)\psi(x;\theta))$, which is not limited to polynomial neurons. Moreover, we present three types of GTNNs: CTNN, QTNN, and TTNN, and they are structurally different from the quadratic networks proposed in previous works. More precisely, compared with the above lines of work on quadratic networks and trainable activation functions, the GTNN framework differs in the following aspects:
\begin{itemize}
    \item \textbf{More general neuron.} The abstract neuron of GTNNs multiplies an affine transformation $Wx+b$ of a conventional neuron by a learnable parametrized function $\psi(x;\theta)$, which can be chosen from a rich family such as polynomials, trigonometric functions, Gaussian functions, etc. This construction allows the network to modulate each neuron's local behavior more flexibly, adapting to the high-frequency or oscillatory nature of the target function, and it also permits a nonlinear dependence on the input $x$. Quadratic networks \citep{Fan2018,Fan2020} appear as special cases of this construction.
    \item \textbf{Approximation properties in Sobolev spaces.} We demonstrate quantitative approximation results of GTNNs for functions in Barron-type function spaces, which include the Sobolev spaces $H^s(\mbR^d)$ with $s = \lfloor d/2 \rfloor + 3$, and the error bound does not depend on the dimension explicitly. It goes beyond the universal approximation or mostly empirical analyses for approximating Barron-type functions in prior quadratic network work.
    \item \textbf{Adaptive parametrized functions.} By choosing specific parametrized functions $\psi(x;\theta)$ (e.g., oscillatory or localized ones), GTNNs can be tailored to approximate special targets arising in high-dimensional learning problems, which is not systematically exploited in existing works, e.g., those on quadratic networks \citep{Fan2018,Fan2020}. We also observe numerically that the proposed CTNNs/QTNNs/TTNNs achieve higher accuracy and better capture local features of the oscillatory target function in least squares regression than conventional DNNs. 
\end{itemize}

\subsection{Organization}
The organization of this paper is as follows. In Section 2, we introduce the general neuron and provide a detailed description of GTNNs. In Section 3, we present the approximation theorems for GTNNs. In Section 4, we present ample numerical experiments, including learning oscillatory functions, the flow of the shock-turbulence problem, multi-peak Gaussian functions and shock waves to validate the advantages of GTNNs. We conclude with some remarks in the final section. 

\section{General Transformation Neural Network} 
\subsection{General neurons}
Abstract neurons are the building blocks of neural networks. Specifically, an abstract neuron is given by
\begin{equation}\label{03}
h(x) = \sigma(w^\top x+b),
\end{equation}
where $x$ is the input vector; $w$ and $b$ are the weight vector and scalar bias, respectively, which are trainable parameters; $\sigma$ is an activation function that applies to its input component-wise and produces an output vector of the same size. 

Clearly, the abstract neuron with ReLU activation $\sigma(t)=\max(0,t)$ is a piecewise linear function, and the domain pieces are separated by the hyperplane $\{x:w^\top x+b=0\}$. In every piece, the ReLU neuron only has a linear shape, which makes it hard to approximate oscillatory functions. Therefore, we expect neurons with more diverse shapes. Our approach is multiplying the affine transformation $w^\top x+b$ by a specific parametrized function $\psi(x;\theta)$, obtaining the general neuron,
\begin{equation}\label{01}
g(x) = \sigma((w^\top x+b)\psi(x;\htheta)).
\end{equation}
Instead of applying to a single affine transformation, the activation of the general neuron applies to a ``general" transformation $(w^\top x+b)\psi(x;\htheta)$. Same as $w$ and $b$, the parameter $\htheta$ of $\psi$ is also trainable in the learning process.

In the case that $\psi$ is smooth, the general neuron $g(x)$ is a piecewise smooth function, and the pieces are separated by surfaces $\{x:(w^\top x+b)\psi(x;\htheta)=0\}$. In every piece where $(w^\top x+b)\psi(x;\htheta)$ takes positive values, $\sigma$ acts as the identity and the general neuron takes $(w^\top x+b)\psi(x;\htheta)$ as its shape function, which is much richer than the linear function $w^\top x+b$. Consequently, the class of general neurons has a larger capacity than the class of abstract neurons.

\subsection{Two-layer architectures}\label{subsec_two_layer_GTNN}
Inspired by the architecture of conventional neural networks, we propose a class of neural network-like functions, i.e., GTNNs, based on the general neurons \eqref{01}. Similar to the two-layer neural networks, which are linear combinations of abstract neurons, the class of two-layer GTNNs is formulated as
\begin{multline}\label{02}
\mathcal{F}_M =\Big\{ \phi:\mbR^{d}\rightarrow\mbR~|~\phi(x; \theta)=\sum_{m=1}^{M} a_m \sigma(( w_m^\top x + b_m)\psi(x; \htheta_m))+c,\\
\forall( a_m, w_m, b_m, \htheta_m )\in\mbR\times\mbR^d\times\mbR\times\Theta,~c\in\mbR\Big\},
\end{multline}
where $M\in\mathbb{N}^+$ is the number of neurons, i.e., the width; $a_m$ and $c$ are the weights and bias of the output layer, respectively; $\Theta$ is the hypothesis space of the parameter $\htheta_m$; $\theta:=\{a_m,w_m,b_m,\htheta_m\}_{m=1}^M\cup\{c\}$ is the set of all free parameters.

When $\sigma$ is the ReLU activation, we may understand the functions $\phi\in\mathcal{F}_M$ by making an analogy with finite and spectral elements. In the case of $\psi\equiv1$, $\phi$ are exactly the two-layer ReLU neural networks, which are piecewise linear continuous functions, i.e., linear finite elements. In this case, the ``mesh" is formed by intersecting the hyperplanes $\{x:w_m^\top x+b_m=0\}_{m=1}^M$, and the ``shape functions" are determined by the components inside $\sigma(\cdot)$, i.e., the linear functions $x\mapsto w_m^\top x+b$. We enrich the class of ReLU neural networks by appending more general $\psi$, leading to the class of two-layer ReLU GTNNs \eqref{02}. For general $\psi$, the ``mesh" still contains edges from the hyperplanes $\{x:w_m^\top x+b_m=0\}_{m=1}^M$ (it also contains edges from the surfaces $\{x:\psi(x;\htheta)=0\}$), and the ``shape functions" are represented by $(w_m^\top x + b_m)\psi(x; \htheta_m)$. For example, if $\psi$ is a $p$-degree polynomial, the ``shape functions" will be ($p+1$)-degree polynomials. So, ReLU GTNNs are like the $hp$ elements, i.e., spectral elements, which are mesh-based functions with a richer class of shape functions.

\subsubsection{Two-layer quadratic transformation neural networks}
Generally, the parametrized function $\psi(x; \htheta)$ can have various forms in different situations. Similar to the polynomial-based spectral elements, we can set $\psi(x;\htheta)$ as polynomials. For example, the general neuron $g(x)$ in \eqref{02} with linear polynomial $\psi$ can be defined by
\begin{equation}\label{26}
g(x) = \sigma(( w^\top x + b )(u^\top x +v)),
\end{equation}
where the vector $u$ and scalar $v$ are both trainable parameters. Note that the activated shape function $( w^\top x + b )(u^\top x +v)$ is a quadratic polynomial. Based on \eqref{26}, we obtain the two-layer quadratic transformation neural network (QTNN)
\begin{equation}\label{08}
\phi(x;\theta) = \sum_{m=1}^{M} a_m \sigma(( w_m^\top x + b_m)(u_m^\top x +v_m))+c.
\end{equation}

\begin{rem}
Previous studies \citep{Fan2020,Fan2018} have already proposed the quadratic neuron: 
\begin{equation}\label{65}
\sigma((w_1^\top x + b_1)(w_2^\top x + b_2) + w_3^\top(x*x) + b_3),
\end{equation}
where $*$ means entry-wise multiplication. And the corresponding two-layer quadratic networks are a linear combination of quadratic neurons, i.e., 
$$f(x) = \sum_{l=1}^{k}\sigma((w_{1,l}^\top x + b_{1,l})(w_{2,l}^\top x + b_{2,l}) + w_{3,l}^\top(x*x) + b_{3,l}).$$ 
Compared with the general neuron \eqref{26}, the quadratic neuron \eqref{65} has an additional quadratic term and a constant term, thus having more parameters. 
\end{rem}

\subsubsection{Two-layer cubic transformation neural networks}
Motivated by the above construction of QTNNs, we further seek a higher-order variant of the parametrized function $\psi(x;\htheta)$ that captures more complex local structures while preserving the architecture of GTNNs \eqref{02}. Specifically, we construct the following two-layer cubic transformation neural network (CTNN):
\begin{equation}\label{07}
 \phi(x;\theta) = \sum_{m=1}^{M} a_m \sigma((w_m^\top x + b_m)(x^\top Q_mx + u_m^\top x +v_m))+c,
\end{equation}
where $Q_m$ is a trainable matrix, $x^\top Q_mx + u_m^\top x +v_m$ is a quadratic form, and the activated function $(w_m^\top x + b_m)(x^\top Q_mx + u_m^\top x +v_m)$ is a cubic polynomial. The neuron of CTNNs keeps the architecture of general neurons \eqref{01}: the affine factor $(w_m^\top x + b_m)$ plays the role of ``mesh" subdivision, while the quadratic form $x^\top Q_mx + u_m^\top x +v_m$ provides a flexible high-order shape function. 

\subsubsection{Two-layer trigonometric transformation neural networks}
Besides taking $\psi(x; \htheta)$ as a polynomial function, we can also set $\psi(x;\theta)$ in the form of a finite Fourier series, i.e., trigonometric polynomials. Then, the general neurons $g(x)$ in \eqref{02} with finite Fourier series $\psi$ can be defined by
\begin{equation}\label{28}
g(x) = \sigma \left(( w^\top x + b )\left(\sum_{n=1}^N c_n\sin(u_n^\top x +v_n) + d_n\cos(u_n^\top x +v_n)\right)\right),
\end{equation}
where the vector $\{u\}_{n=1}^N$, scalar $\{v\}_{n=1}^N$, weight $c_n$ and $d_n$ are all trainable parameters. In practical applications, to avoid having too many parameters, we can just set $\psi(x;\theta)$ as $\sin(u^\top x +v)$ or $\cos(u^\top x +v)$. Based on \eqref{28}, we can present the two-layer trigonometric transformation neural network (TTNN) architecture: 
\begin{equation}\label{29}
\phi(x;\theta) = \sum_{m=1}^{M} a_m \sigma\left(( w_m^\top x + b_m)\left(\sum_{n=1}^N c^m_n\sin(u_{m,n}^\top x +v_{m,n}) + d^m_n\cos(u_{m,n}^\top x +v_{m,n})\right)\right)+c.
\end{equation}

More choices of $\psi(x;\htheta)$ include wavelet functions, Gaussian functions, etc., which are not discussed in detail in this paper.

\subsection{Deep architectures}\label{subsec_deep_GTNN}
Moreover, we propose deep GTNNs, referring to the architecture of (deep) fully connected neural networks (FNNs). Letting $L\geq2$ be the depth and $M_\ell \in \mathbb{N}^+$ ($\ell = 1, \ldots, L-1$) be the widths of hidden layers, the FNN is the composition of $L-1$ layers of abstract neurons; namely,
\begin{equation}\label{04}
\phi_{\mathrm{FNN}}(x; \theta) = a^\top h_{L-1} \circ h_{L-2} \circ \cdots \circ h_1(x)+c \quad \text{for } x \in \mbR^d,
\end{equation}
where $a\in\mbR^{M_{L-1}}$ and $c\in\mbR$ are the weight vector and bias scalar of the output layer, respectively; $h_\ell(x_\ell)=\sigma(W_\ell x_\ell + b_\ell)$ with $W_\ell \in \mbR^{M_\ell \times M_{\ell-1}}$ ($M_0:= d$) and $b_\ell \in \mbR^{M_\ell}$ is the vector-valued abstract neuron in the $\ell$-th layer; $\theta:=\{W_\ell, b_\ell\}_{\ell=1}^{L-1}\cup\{a,c\}$ is the set of all trainable parameters. 

\begin{figure}[htbp]
\centering
\includegraphics[scale=0.5,trim=120 145 140 90,clip]{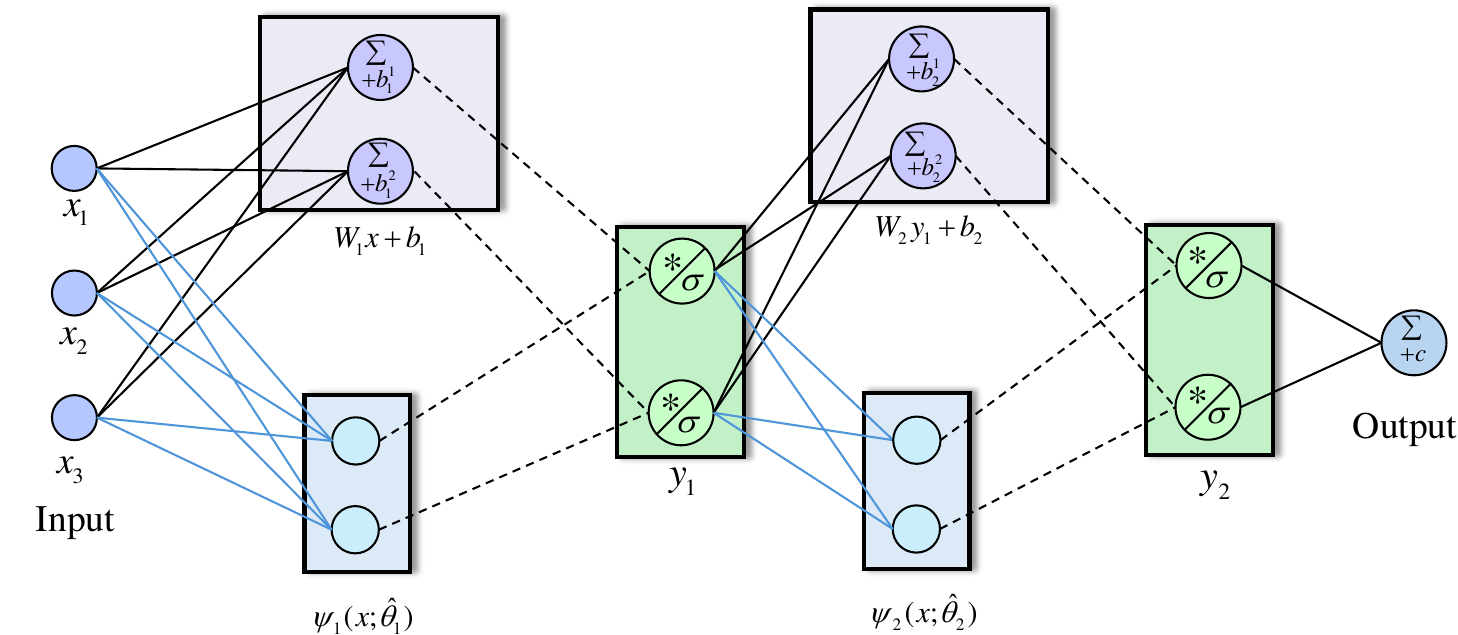}
\caption{\em The architecture of a three-layer GTNN with $L=3$ and $M_1=M_2=2$ (the black dashed lines mean the scalar multiplication, and the blue solid lines mean general operations implemented by $\psi$).}
\label{fig:three-layer GTNN}
\end{figure}

Similarly, we can define the vector-valued general neuron as
\begin{equation}
g(x)= \sigma((Wx+b)\ast\psi(x; \htheta)),
\end{equation}
where $\ast$ means entry-wise multiplication; $W$ and $b$ are the weight matrix and bias vector, respectively. Then, we define the deep GTNN as 
\begin{equation}\label{09}
\phi(x; \theta) := a^\top g_{L-1} \circ g_{L-2} \circ \cdots \circ g_1(x)+c \quad \text{for } x \in \mbR^d,
\end{equation}
where $g_\ell(x_\ell)= \sigma((W_\ell x_\ell+b_\ell)\ast\psi_\ell(x_\ell; \htheta_\ell))$ is a general neuron for every $\ell$. In Figure \ref{fig:three-layer GTNN}, we present the architecture of the GTNN with $L=3$.

\subsubsection{Deep quadratic transformation neural networks}
\begin{figure}[htbp]
\centering
\includegraphics[scale=0.5,trim=120 100 140 90,clip]{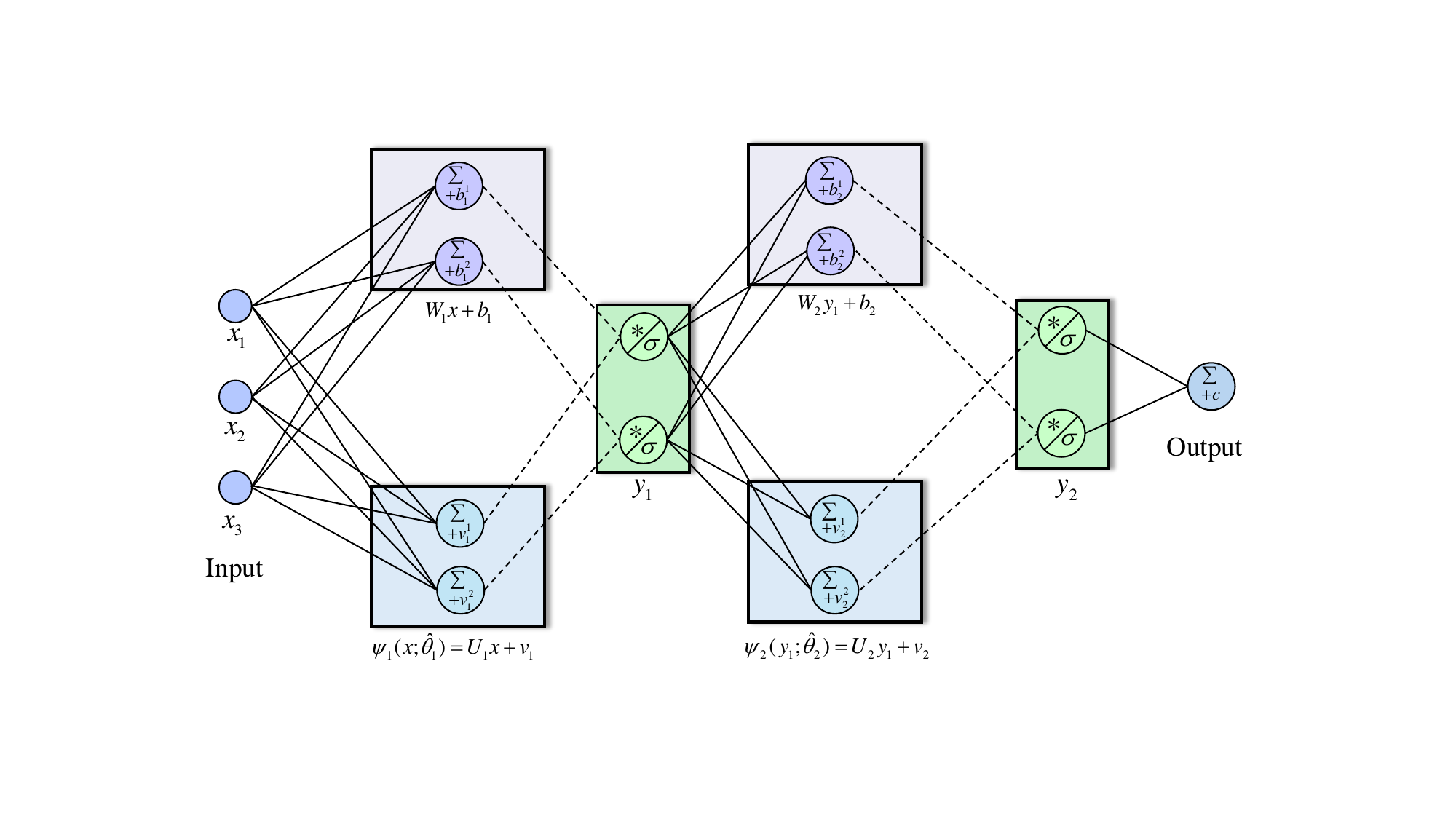}
\caption{\em The architecture of a three-layer QTNN with $L=3$ and $M_1=M_2=2$ (the black dashed lines mean the scalar multiplication).}
\label{fig:three-layer QTNN}
\end{figure}

A simple example of GTNNs is to take $\psi_\ell$ as a linear polynomial, i.e.,
\begin{equation}\label{23}
\psi_\ell(x; \htheta_\ell) = U_\ell x +v_\ell,
\end{equation}
where $\htheta_\ell=\{U_\ell,v_\ell\}$ is trainable. In this case, $(W_\ell x_\ell+b_\ell)\ast\psi(x_\ell; \htheta_\ell)$ is a vector-valued quadratic polynomial. Therefore, the deep GTNN \eqref{09} with $\psi_\ell$ given by \eqref{23} is named as the deep QTNN. In Figure \ref{fig:three-layer QTNN}, we present the architecture of the QTNN with $L=3$ and $M_1=M_2=2$. Besides, we note that when $L=2$, the deep QTNN is exactly the two-layer QTNN \eqref{08}.

\subsubsection{Deep cubic transformation neural networks}
Similarly, one can set $\psi_\ell$ as a reduced quadratic polynomial
\begin{equation}\label{24}
\psi_\ell(x; \htheta_\ell) = [x^\top Q_\ell^i x]_{M_\ell\times1} +U_\ell x +v_\ell,
\end{equation}
where $Q_\ell^i\in\mbR^{M_{\ell-1}\times M_{\ell-1}}$ ($i=1,\dots,M_\ell$), $U_\ell\in\mbR^{M_\ell \times M_{\ell-1}}$, $v_\ell \in \mbR^{M_\ell}$, and $\htheta_\ell=\{Q_\ell^i,U_\ell,v_\ell\}$ is the set of trainable parameters. In this case, $(W_\ell x_\ell+b_\ell)\ast\psi(x_\ell; \htheta_\ell)$ is a vector-valued cubic polynomial, and we name the deep GTNN \eqref{09} with $\psi_\ell$ \eqref{24} as the deep CTNN. When $L=2$, the deep CTNN is equivalent to the two-layer CTNN \eqref{07}.

\subsubsection{Deep trigonometric transformation neural networks}
Moreover, $\psi_\ell$ can also be a finite Fourier series
\begin{equation}\label{52}
\psi_\ell(x; \htheta_\ell) = \sum_{n=1}^N C^\ell_n \ast\sin(U_n^\ell x +v_n^\ell) + D^\ell_n\ast\cos(U_n^\ell x +v_n^\ell),
\end{equation}
where $U_\ell\in\mbR^{M_\ell \times M_{\ell-1}}$, $v_\ell, C^\ell_n, D^\ell_n \in \mbR^{M_\ell}$, and $\htheta_\ell=\{C^\ell_n, D^\ell_n, U^\ell_n,v^\ell_n\}$ is the set of trainable parameters. In this case, $(W_\ell x_\ell+b_\ell)\ast\psi(x_\ell; \htheta_\ell)$ is a vector-valued finite Fourier series multiplied by an affine transformation, and we name the deep GTNN \eqref{09} with $\psi_\ell$ \eqref{52} as the deep TTNN. When $L=2$, the deep TTNN is equivalent to the two-layer TTNN \eqref{29}.

Generally, the form of $\psi_\ell$ can be various for different hidden layers in \eqref{52}. Although the class of GTNNs has a richer capacity than the class of DNNs with the same size, the former has more free parameters to be tuned. For example, the abstract neuron $h(x)=\sigma(w^\top x+b)$ from $\mbR^d$ to $\mbR$ has parameters $w\in\mbR^d$ and $b\in\mbR$, so the total number is $d+1$. Comparatively, the QTNN neuron in \eqref{08} has $2d+2$ parameters, and the CTNN neuron in \eqref{07} has $d^2+2d+2$ parameters. Therefore, GTNNs will require more storage than DNNs of the same size in practical computations. However, compared with conventional DNNs, GTNNs can achieve much higher accuracy in learning problems with the same number of trainable parameters (see Section \ref{sec_numerical_examples}).

We remark that the difference between finite/spectral elements and ReLU DNNs/GTNNs lies in the variability of the mesh. For the former, one must build the mesh before setting up the parametrized form of functions, and the mesh is typically fixed in practical computations. But for ReLU DNNs/GTNNs, the ``mesh" is generated by the intersection of the hyperplanes $\{x:Wx+b=0\}$, where $W$ and $b$ are trainable parameters; therefore, the ``mesh" is varying and can be optimized in the learning process. Moreover, compared with ReLU DNNs whose ``shape functions" are linear functions, GTNNs have a richer class of ``shape functions", so GTNNs usually have larger capacities than DNNs with the same depth and width, and are more capable of approximating high-frequency or oscillatory functions. 

\begin{figure}[htbp]
\centering
\subfloat[Piecewise quadratic polynomial]{
\includegraphics[scale=0.35]{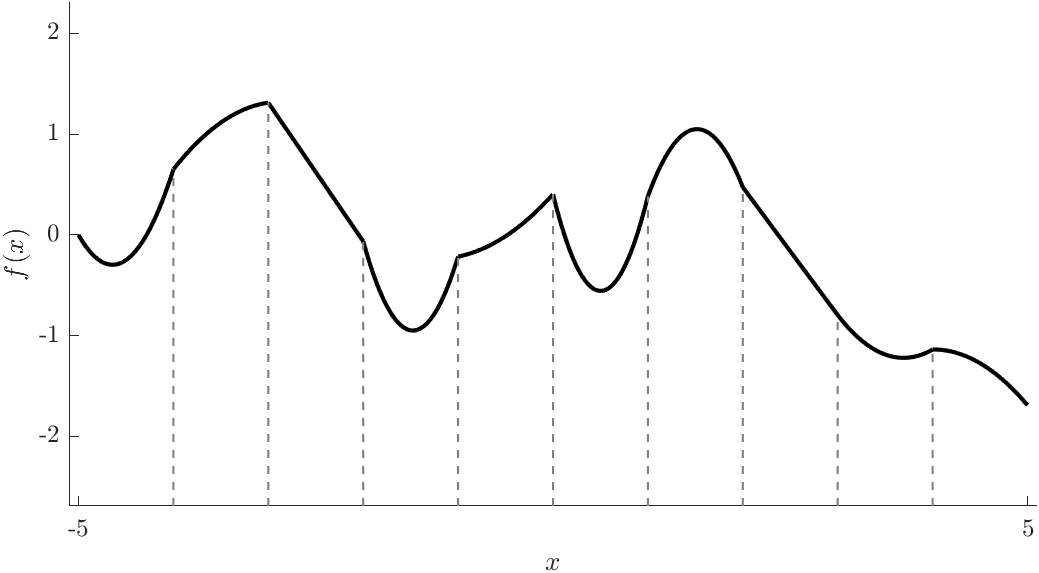}}\hspace{10pt}
\subfloat[Two-layer QTNN]{
\includegraphics[scale=0.35]{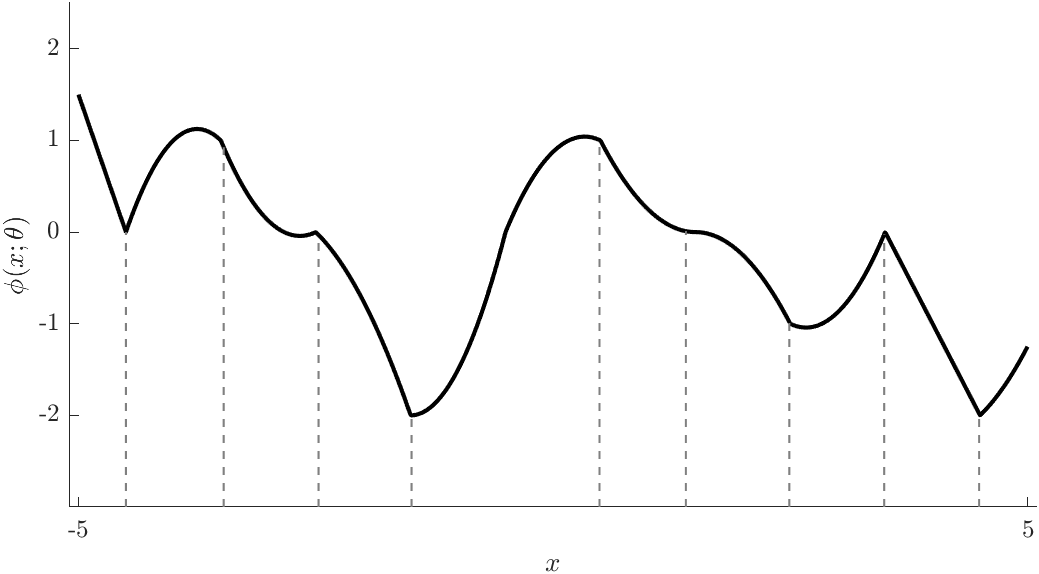}}\\
\subfloat[Piecewise cubic polynomial]{
\includegraphics[scale=0.35]{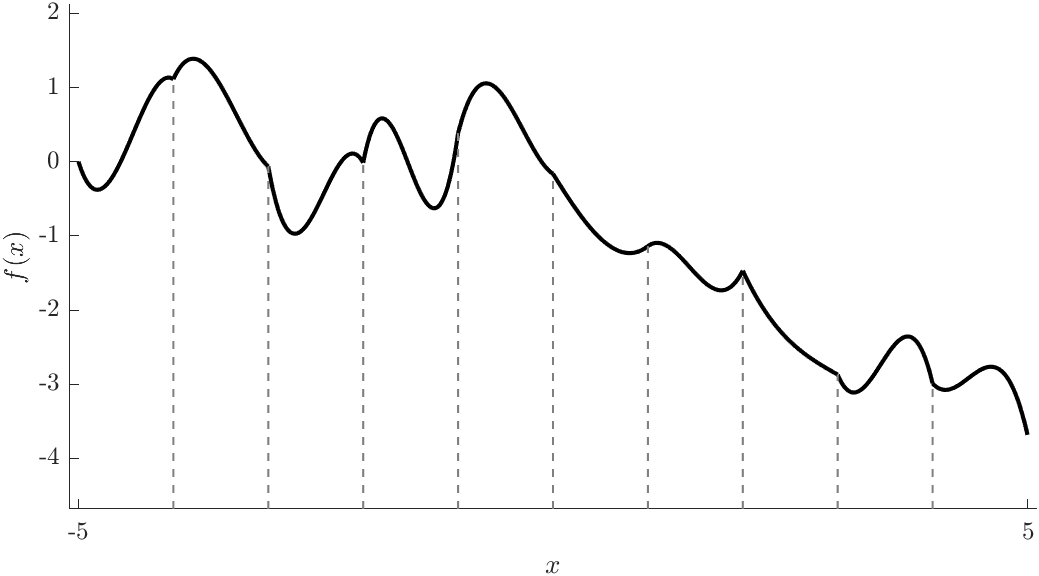}}\hspace{10pt}
\subfloat[Two-layer CTNN]{
\includegraphics[scale=0.35]{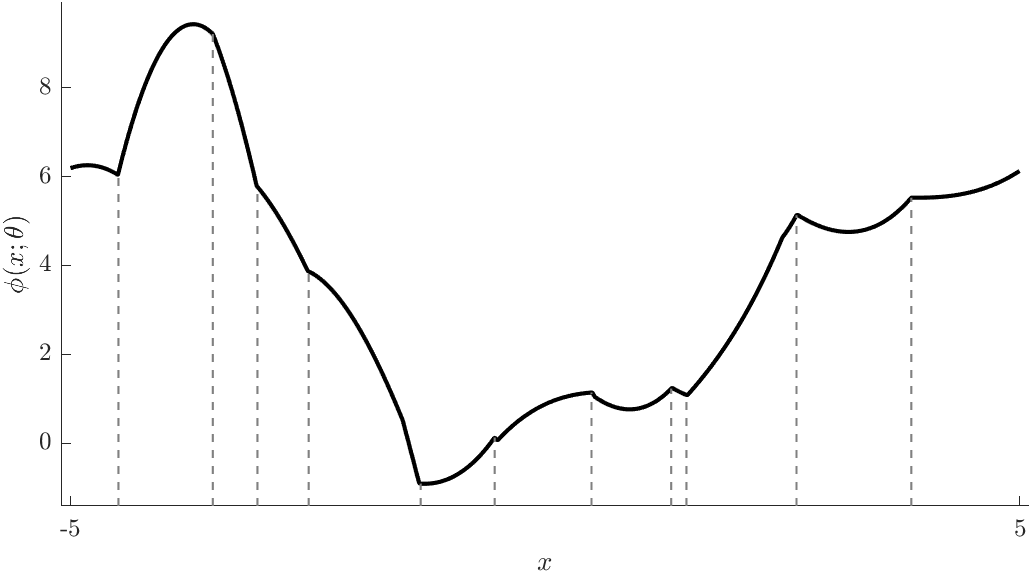}}
\caption{\em Graphs of the piecewise quadratic/cubic polynomial (with fixed equispaced grids) and the two-layer QTNN/CTNN (with adaptive trainable grids).}
\label{Fig_01}
\end{figure}

In Figure \ref{Fig_01}, we present simple examples of 1-D piecewise quadratic and cubic polynomials on equispaced grids, as well as the two-layer QTNNs and CTNNs on the interval $[-5,5]$. The coefficients of polynomials and the parameters of the QTNN/CTNN are randomly set from $U(-2,2)$. Clearly, the ``grids" of the QTNN and CTNN are randomly generated, not equispaced, and trainable.

\section{Approximation Property} 
This section presents a universal approximation theorem and the approximation rate estimation for GTNNs. Our analysis is based on the approximation property of DNNs for continuous functions and Barron-type functions, and deep DNNs. We will prove that the two-layer GTNN can approximate continuous functions with any desired degree of accuracy. We also derive error bounds for approximating Barron-type functions and error bounds of deep architectures. These results not only validate the approximation property of GTNNs but also provide insights into their practical applicability.  

\subsection{Universal approximation theorems}\label{UA}
\subsubsection{General approximation}
First, the universal approximation theorem for the two-layer GTNN is given as follows. 
\begin{thm}\label{thm1}
Suppose $\sigma\in C(\mbR)$ is not a polynomial and $\psi(\cdot;\htheta_0)$ is a non-zero constant for some $\htheta_0\in\Theta$. Let $d\geq1$ and $\Omega \subset \mbR^d$ be a compact set. For any $f \in C(\Omega)$ and $\epsilon>0$, there exists $M\in\mathbb{Z}^+$ and $a,b \in\mbR^{M}$, $W \in\mbR^{M \times d }$, $c \in \mbR$ such that the GTNN $\phi(x) = a^\top \sigma((Wx + b) \ast \psi(x;\htheta))  +c$ satisfies
\begin{equation}
\left|f(x) - \phi(x)\right| < \epsilon,
\end{equation}
for all $x \in \Omega$.
\end{thm}
\begin{proof}
Taking $\htheta=\htheta_0$ in $\phi(x)$, then $\phi(x)=a^\top \sigma(Wx + b) + c$ is simplified to a two-layer FNN with activation $\sigma$, which belongs to the two-layer FNN class $\mathcal{M}(\sigma)=\left\{ \phi(x)=a^\top \sigma(Wx + b) + c:~ a, b \in\mbR^{M}, W \in\mbR^{M \times d }, c \in \mbR \right\}$. We now recall the universal approximation theorem \citep[Theorem~3.1]{Pinkus1999}: if $\sigma\in C(\mathbb{R})$ is not a polynomial, then $\mathcal{M}(\sigma)$ is dense in $C(\mathbb{R}^d)$ in the topology of uniform convergence on compact sets. Therefore, for any compact $\Omega\subset\mathbb{R}^d$, any $f\in C(\Omega)$ and any $\epsilon>0$, there exists $\phi\in \mathcal{M}(\sigma)$ such that
\begin{equation}
\sup_{x\in\Omega}|f(x)-\phi(x)|<\epsilon,
\end{equation}
which proves the claim.
\end{proof}

\subsubsection{Application to CTNNs, QTNNs and TTNNs}
Now we consider three types of GTNN: QTNN, CTNN, and TTNN. For QTNNs $\phi(x) = a^\top \sigma((Wx + b) \ast (Ux +v))  +c$, we can take $U$ as the zero matrix and $v = [1~\dots~1]^\top$. The next proof is similar to Theorem \ref{thm1}. For CTNNs $\phi(x) = a^\top \sigma((Wx + b) \ast (Q(x)+Ux+v)) +c$, where $Q(x) = [x^\top V_1 x~\dots~x^\top V_M x]^\top$ and $V_i\in\mbR^{d \times d }$, we let  $V_i$ ($i=1,\dots,M$) as the zero matrices and the same conclusion can be drawn for CTNNs. The same results can be obtained for TTNN through discussions similar to CTNNs, because one can choose parameters so that $\psi(x;\hat\theta)$ becomes a nonzero constant (e.g., using $\cos(0)=1$).

If $\sigma$ is a polynomial, Theorem \ref{thm1} may not hold. For example, if $\sigma$ is a polynomial of degree $m$, then, for every choice of $W$, $U\in\mbR^{M \times d }$, $b$, $v\in\mbR^M$ of QTNNs, $\phi(x)$ is a multivariate polynomial of total degree at most $2m$. But the space of polynomials of degree at most $2m$ does not span $C(\mbR)$. The above theorems indicate that if $\psi(x)$ and $\sigma$ are chosen appropriately, the GTNN is a universal approximation to continuous functions with any desired degree of accuracy. 

\subsection{Error bounds for Barron-type functions} 
In this section, we derive error bounds for two-layer GTNNs for Barron-type target functions and conduct a detailed analysis of two specific GTNN architectures (CTNN and QTNN), deriving error bounds of their approximation. Roughly, Barron functions can be analogous to infinitely wide two-layer neural networks. The bounds do not depend on the dimension explicitly, so the approximation lessens or overcomes the curse of dimensionality in some sense. 

\subsubsection{General approximation}\label{General}
In this part, we first consider the approximation of GTNNs. Following \eqref{02}, the class of the two-layer ReLU GTNN of width $M$ can be reformulated as
\begin{multline}
\mathcal{F}_{M,G} =\Big\{ \phi:\mbR^{d}\rightarrow\mbR~|~\phi(x; \theta)=\frac{1}{M}\sum_{m=1}^{M} a_m \sigma(( w_m^\top x + b_m)\psi(x; \htheta_m)),\\
\forall( a_m, w_m, b_m, \htheta_m )\in\mathcal{D}:=\mbR\times\mbR^d\times\mbR\times\Theta\Big\},
\end{multline}
where $\sigma$ is the ReLU activation. Here, we omit the bias of the output layer. For a function $f: \mbR^{d} \to \mbR$, we say that a probability measure $\pi$ on $\mathcal{D}$ represents $f$ if 
\begin{multline}\label{64}
f(x) = f_{\pi}(x) = \mathbb{E}_{\pi}\left[a\sigma\left((w^\top x + b)\psi(x; \htheta)\right)\right]  \\
=\int_{\mathcal{D}} a\sigma\left((w^\top x + b)\psi(x; \htheta)\right) \pi(\td a, \td w, \td b, \td\htheta), ~ \forall x \in \mbR^{d},
\end{multline}
whenever the integral is well defined. We set $\Pi_f$ to denote the set of all probability measures $\pi$ such that $f_{\pi}(x) = f(x)$ almost everywhere. For any $\psi(x;\htheta)$ which satisfies $M_\psi(\htheta):=\max(\sup_{x\in \Omega}|\psi(x;\htheta)|, \sup_{x\in \Omega}|\nabla\psi(x;\htheta)|)<\infty$, we can define the following Barron-type constant 
\begin{multline}\label{B_norm_GTNN}
\mathcal{C}_{\psi}(f) := 
\left(\inf_{\pi\in\Pi_{f}} \int_{\mathcal{D}} a^{2}M_\psi(\htheta)^2(|w| + |b|)^{2}\pi(\td a, \td w, \td b, \td\htheta) \right)^{1/2} \\
=\left(\inf_{\pi \in \Pi_{f}} \mathbb{E}_{\pi}\left[a^2M_\psi(\htheta)^2\left(|w| + |b|\right)^2\right]\right)^{1/2},
\end{multline}
where $|w|=(w\cdot w)^{\frac{1}{2}}$ and the infimum of the empty set is considered to be $+\infty$. The class of all functions $f$ with $\mathcal{C}_{\psi}(f)<\infty$ is denoted by $\mathcal{B}$, which can be viewed as a Barron-type space \citep{E2020,E2022,E2019}.  Recall that the $H^1$ Sobolev norm of a function $f$ over a domain $\Omega$ is defined by 
\begin{equation}
\|f\|^2_{H^{1}(\Omega)} :=  \int_{\Omega} |f|^2 + |\nabla f|^2 \td x.
\end{equation}
Then we have the following result.
\begin{thm}\label{thm2}
Let $\Omega \subset \mbR^{d}$ be a bounded set. Suppose that $\psi(x;\htheta)\in C^1(\Omega)$. Then for any $f\in \mathcal{B}$, there exists some $\phi \in \mathcal{F}_{M,G}$ such that  
\begin{equation}
\|\phi - f\|_{H^{1}(\Omega)} \leq C_\Omega \mathcal{C}_{\psi}(f) M^{-\frac{1}{2}},
\end{equation}
where $C_{\Omega} := \sqrt{5}|\Omega|^{\frac{1}{2}}\max(\sup_{x\in \Omega}|x|, 1)$. 
\end{thm}
\begin{proof}
Since $f\in \mathcal{B}$, we have $\mathcal{C}_{\psi}(f) <\infty$, and there exists some probability measure $\pi$ such that $f_\pi = f$ almost everywhere and  $\mathbb{E}_{\pi}\left[a^2\left(|w| + |b|\right)^2 M_\psi^2 \right] \leq 2\mathcal{C}^2_{\psi}(f)$. For $(a,w,b,\htheta)\in \mathcal{D}$, using the fact $|\sigma(t)| \leq |t|$ and $|\sigma'(t)| = \chi_{t\geq0}$, we have
\begin{multline}\label{10}
\|a\sigma((w^\top x + b)\psi(x;\htheta))\|^2_{H^{1}(\Omega)}=\int_{\Omega} |a\sigma((w^\top x + b)\psi(x;\htheta))|^2 + |\nabla (a\sigma((w^\top x + b)\psi(x;\htheta))|^2 \td x \\
\leq \int_{\Omega} a^2(w^\top x + b)^2\psi^2(x;\htheta) +a^2 \chi_{t\geq0} |\nabla((w^\top x + b)\psi^2(x;\htheta))|^2 \td x \\
\leq \int_{\Omega} a^2\left(|w||x| + |b|\right)^2(\sup_{x\in\Omega}|\psi(x;\htheta)|)^2 +a^2(\sup_{x\in\Omega}|\psi(x;\htheta)||w|+\sup_{x\in\Omega}|\nabla \psi(x;\htheta)||w^\top x+b|)^2  \td x\\
\leq \int_{\Omega} a^2\left(|w||x| + |b|\right)^2(\sup_{x\in\Omega}|\psi(x;\htheta)|)^2 +2a^2((\sup_{x\in\Omega}|\psi(x;\htheta)|)^2|w|^2+(\sup_{x\in\Omega}|\nabla \psi(x;\htheta)|)^2|(|w| |x|+|b|)^2)  \td x. 
\end{multline} 
Since $\Omega$ is a bounded set,  for any $x \in \Omega $ we have $|x| \leq R_{\Omega} $, where $R_{\Omega} := \max\left(\sup_{x\in \Omega}|x|, 1\right)$, and since $\psi(x;\htheta)\in C^1(\Omega)$, we have $M_\psi=\max(\sup_{x\in \Omega}|\psi(x;\htheta)|, \sup_{x\in \Omega}|\nabla\psi(x;\htheta)|) <\infty$ . Then we have
\begin{multline}\label{11}
\int_{\Omega} a^2\left(|w||x| + |b|\right)^2(\sup_{x\in\Omega}|\psi(x;\htheta)|)^2 +2a^2((\sup_{x\in\Omega}|\psi(x;\htheta)|)^2|w|^2+(\sup_{x\in\Omega}|\nabla \psi(x;\htheta)|)^2|(|w||x|+|b|)^2)  \td x \\
\leq \int_{\Omega} a^2\left(|w|R_{\Omega} + |b|\right)^2(\sup_{x\in\Omega}|\psi(x;\htheta)|)^2 +2a^2((\sup_{x\in\Omega}|\psi(x;\htheta)|)^2|w|^2+(\sup_{x\in\Omega}|\nabla \psi(x;\htheta)|)^2|(|w| R_{\Omega}+|b|)^2) \td x  \\
\leq \int_{\Omega} 5M_\psi^2 a^2 \left(|w|R_{\Omega} + |b|\right)^2 \td x 
\leq \int_{\Omega} 5R_{\Omega}^2M_\psi^2 a^2\left(|w| + |b|\right)^2 \td x  = 5R_{\Omega}^2|\Omega| M_\psi^2 a^2\left(|w| + |b|\right)^2.
\end{multline}
Combining \eqref{10} and \eqref{11}, we have
\begin{equation}\label{12}
\mathbb{E}_{\pi}\left[\|a\sigma((w^\top x + b)\psi(x;\htheta))\|_{H^{1}(\Omega)}^2\right] 
\leq 5R_{\Omega}^2 |\Omega|\mathbb{E}_{\pi}\left[a^2 M_\psi^2 \left(|w| + |b|\right)^2 \right]
\leq 5R_{\Omega}^2 |\Omega|\mathcal{C}^2_{\mathcal{\psi}}(f). 
\end{equation}

On the other hand, the mapping
\begin{equation}
    \mathcal{D}\to\mbR,\quad (a,w,b,\htheta)\to a\sigma((w^\top x + b) \psi(x;\htheta))
\end{equation}
is continuous and hence Bochner measurable. Also, by Cauchy-Schwartz inequality, (\ref{12}) leads to
\begin{equation}
    \mathbb{E}_{\pi}\left[\|a\sigma((w^\top x + b)\psi(x;\htheta))\|_{H^{1}(\Omega)}\right] \leq \left(\mathbb{E}_{\pi}\left[\|a\sigma((w^\top x + b)\psi(x;\htheta))\|_{H^{1}(\Omega)}^2\right]\right)^{1/2},
\end{equation}
which implies the integral $\int_{\mathcal{D}} a\sigma\left((w^\top x + b)\psi(x;\htheta)\right) \pi(\td a, \td w, \td b, \td \htheta)$ is a Bochner integral.

We note the fact that if $\xi_{1},\ldots,\xi_{M}$ are i.i.d random variables with a random variable $\xi$, then
\begin{multline*}
\quad \mathbb{E}\left(\frac{1}{M}\sum_{i=1}^{M}\xi_{i} - \mathbb{E}\xi\right)^{2} 
      = \mathbb{E}\left(\frac{1}{M}\sum_{i=1}^{M}\left(\xi_{i} - \mathbb{E}\xi\right)\right)^{2} \\
= \frac{1}{M^{2}}\left(\sum_{i=1}^{M}\mathbb{E}\left(\xi_{i} - \mathbb{E}\xi\right)^{2} 
      + \sum_{1\leq i<j\leq M}\mathbb{E}\left(\xi_{i} - \mathbb{E}\xi\right)\cdot\mathbb{E}\left(\xi_{j} - \mathbb{E}\xi\right)\right) \\
    = \frac{1}{M}\sum_{i=1}^{M}\mathbb{E}\left(\xi - \mathbb{E}\xi\right)^{2} 
= \frac{1}{M}\mathbb{E}\xi^{2} - \frac{1}{M}\left(\mathbb{E}\xi\right)^{2} 
      \leq \frac{1}{M}\mathbb{E}\xi^{2}.
\end{multline*}
By a similar argument, for i.i.d random variables $\{(a_{i}, w_{i}, b_{i}, u_{i}, v_{i})\}$ from $\pi$, we have
\begin{equation}
    \mathbb{E}_{\pi^M}\left[\left\|\frac{1}{M} \sum_{i=1}^{M} a_{i} \sigma(( w_{i}^\top  x + b_{i} )\psi(x;\htheta_i)) -f_{\pi}(x)\right\|^2_{H^{1}(\Omega)}\right]
    \leq \frac{1}{M}\mathbb{E}_{\pi}\left[\left\|a\sigma((w^\top x + b)\psi(x;\htheta))\right\|^2_{H^{1}(\Omega)}\right]
\end{equation}
In particular, there exists a set of samples $\{(a_{i}, w_{i}, b_{i}, \htheta_{i})\}$ such that 
\begin{equation}\label{13}
    \left\|\frac{1}{M} \sum_{i=1}^{M} a_{i} \sigma(( w_{i}^\top  x + b_{i} )\psi(x;\htheta_i)) -f_{\pi}(x)\right \|^2_{H^{1}(\Omega)} 
    \leq \frac{1}{M}\mathbb{E}_{\pi}\left[\left\|a\sigma((w^\top x + b)\psi(x;\htheta))\right\|^2_{H^{1}(\Omega)}\right]
\end{equation}
Let $\phi = \frac{1}{M} \sum_{i=1}^{M} a_{i} \sigma(( w_{i}^\top  x + b_{i} )\psi(x;\htheta_i))$. Combining (\ref{13}) and (\ref{12}), we obtain 
\begin{equation}
    \|\phi - f\|_{H^{1}(\Omega)} \leq C_\Omega\mathcal{C}_{\mathcal{\psi}}(f) M^{-\frac{1}{2}}.
\end{equation}
\end{proof}
This theorem converts the qualitative integral representation into a quantitative approximation bound for GTNNs of width $M$. It implies that a class of functions with finite Barron-type constant can be efficiently approximated by two-layer GTNNs in the $H^1$-sense. In what follows, we will present some results to illustrate that the Sobolev space $H^s(\mbR^d)=\{f:\mbR^d\rightarrow\mbR~|~ D^\alpha f\in L^2(\mbR^d),~ \forall \|\alpha\|_1\leq s\}$ is contained in $\mathcal{B}$, where $\alpha=(\alpha_1,\dots,\alpha_d)$ is a multi-index with $D^\alpha f=\partial^{\alpha_1}_{x_1}\cdots\partial^{\alpha_d}_{x_d} f$ and $\|\alpha\|_1=\alpha_1+\dots+\alpha_d$.

Let $\Omega \subset \mbR^{d}$ be a bounded set. The pioneering works \citep{Barron1993,Klusowski2018} have already proposed the Barron functions $f:\Omega\rightarrow\mbR$, which have a finite Fourier–moment functional 
\begin{equation}\label{18}
   \Gamma(f) := \inf_{\hat{f}}\int_{\mbR^d} \|\omega\|_2^2|\hat{f}(\omega)|\td \omega <\infty,
\end{equation}
where $\hat{f}(\omega)$ is the Fourier transform of an extension of $f$ from $\Omega$ to $\mbR^d$ and the infimum is taken over the Fourier transform of $f$ when $f\in C(\mbR^d)$. We introduce $\Gamma(\cdot)$ as an auxiliary functional that will serve as a key tool in the proof of the approximation properties below. Next, we will prove that sufficiently smooth functions in the Sobolev space  $H^s(\mbR^d)$ satisfy \eqref{18}. 
\begin{lem} \label{lem1}
Let $d \ge 1$ and set $s = \lfloor d/2 \rfloor + 3$. For any function $f\in C(\mbR^d)$ with $D^\alpha f\in L^2(\mbR^d)$ for all $\|\alpha\|_1=2$ and $\|\alpha\|_1=s$, there exists a constant $C(d,s) > 0$ such that
\begin{equation}
\Gamma(f) \leq C(d,s)\left(\sum\limits_{\|\alpha\|=2}\|D^\alpha f\|_{L^2}^2 + \sum\limits_{\|\beta\|_1=s}\|D^\beta f\|_{L^2}^2\right)^{\frac{1}{2}},
\end{equation}
where $C(d,s)= (2\pi)^{d/2} \left( \int_{\mbR^d}\frac{\td \omega}{1+\|\omega\|_2^{2s-4}}\right)^{\frac{1}{2}} \max\limits_{|\beta|=s}\binom{s}{\beta}$. 
\end{lem}

The proof of Lemma \ref{lem1} can be found in the Appendix \ref{prflem1}. Obviously, Lemma \ref{lem1} shows that the functions in Sobolev space $H^s(\mbR^d)$ with $s = \lfloor d/2 \rfloor + 3$ satisfy \eqref{18}. The choice $s = \lfloor d/2 \rfloor + 3$ ensures the integrability condition $2s-4 > d$, which is crucial for the convergence of the integral term. From the perspective of high-dimensional approximation, this provides a regularity threshold under which Sobolev functions can be treated as Barron functions, offering guidance for the design and analysis of neural network models. To illustrate that $\mathcal{B}$ contains a large number of functions, we will introduce how $H^s(\mbR^d)$ can be represented by functions in $\mathcal{B}$.

\begin{lem}\label{lem2}
Let $\Omega \subset \mbR^d$ be a bounded set. Assume that $\psi(\cdot;\htheta)$ is nonnegative and continuous for all $\htheta\in\Theta$, and there exists a probability measure $\nu$ on $\Theta$ which satisfies the partition of unity
\begin{equation}\label{51}
\int_{\Theta} \psi(x;\htheta) \nu(\td \htheta) =1,~ \forall x\in \Omega.
\end{equation}
Then, for any $f\in H^s(\mbR^d)$ with $s = \lfloor d/2 \rfloor + 3$, there exists a finite signed measure $\rho$ on $\mathcal{D}_0 = \mbR^d\times\mbR\times\Theta$ such that
\begin{equation}
f(x) = \int_{\mathcal{D}_0} \sigma((w^\top x+b)\psi(x;\htheta)) \rho(\td w, \td b, \td \htheta),~ \forall x\in \Omega,
\end{equation}
where $\sigma$ is the ReLU activation.
\end{lem}

The proof of Lemma \ref{lem2} can be found in the Appendix \ref{prflem2}. This result shows that Sobolev functions $f\in H^s(\mbR^d)$ with $s = \lfloor d/2 \rfloor + 3$ can be represented exactly in a GTNN-type integral form on a bounded set, in terms of a finite signed measure on the parameter domain. In other words, at the level of infinite width, the abstract GTNN architecture is rich enough to represent a large and classical function class, such as $H^s(\mbR^d)$, once an appropriate signed measure on the parameter domain is allowed. To further illustrate the relationship between $H^s(\mbR^d)$ and Barron-type space $\mathcal{B}$, we have the following result

\begin{thm}\label{thm3}
Let $d\geq1$ and $\Omega\subset\mbR^d$ be a bounded set. Assume that $\psi(\cdot;\htheta)$ is nonnegative and continuous for all $\htheta\in\Theta$, satisfying the partition of unity \eqref{51} and 
\begin{align}\label{53}
J_{\psi,\nu}:=\int_{\Theta} M_\psi(\htheta)^2\,\nu(\td \htheta) < \infty
\end{align}
for some probability measure $\nu$ on $\Theta$. Then $H^s(\mbR^d)\subset\mathcal{B}$, where $s = \lfloor d/2 \rfloor + 3$.
\end{thm}

The proof of Theorem \ref{thm3} can be found in the Appendix \ref{prfthm3}. Combining Theorem \ref{thm2} with Theorem \ref{thm3}, we conclude that any $f \in H^s(\mbR^d)$ with $s = \lfloor d/2 \rfloor + 3$ belongs to $\mathcal{B}$ and hence it can be efficiently approximated by two-layer GTNNs in the $H^1$-sense. From a practical point of view, Theorem \ref{thm3} provides theoretical guidance for choosing the network width in high-dimensional approximation and PDE-solving problems, as they link the required model size to intrinsic regularity and Barron-type properties of the underlying target function. 

\subsubsection{Application to QTNNs}
In what follows, we verify that the parametrized functions $\psi$ of QTNNs indeed satisfy the assumptions of Theorem \ref{thm3}. More specifically, we construct a suitable parameter set $\Theta$ and probability measure for the parametrized function $\psi$ of QTNNs so that it satisfies the assumptions on $\Omega$ required by Theorem \ref{thm3}. 

Recall that the parametrized function $\psi$ of QTNNs is $\psi(x;\htheta)=u^\top x+v$. We first consider the simple case in which the hypothesis space of parameters $\Theta$ consists of only two elements. Since $\Omega$ is bounded, we can take a vector $u_0\in\mbR^d$ such that $|u_0|\le \frac{1}{\sqrt{d}\max_{x\in\Omega}\|x\|_\infty}$, and set $\htheta_1=(u_0,1)\in\Theta$, $\htheta_2=(-u_0,1)\in\Theta$. We define the probability measure $\nu:=\tfrac12\delta_{\htheta_1}+\tfrac12\delta_{\htheta_2}$, where $\delta_{\htheta}$ is the Dirac mass at $\htheta$. We take $\mathcal{D} = \mbR^{d} \times \mbR \times \Theta $ where $\Theta = \{\htheta_1,\htheta_2\}$.
Then for all $x\in\Omega$, we have
\begin{align}
& \psi(x;\htheta_1)=u^\top_0 x+1\ge 1-\sqrt{d}\max_{x\in\Omega}\|x\|_\infty|u_0|\ge 0,\\
& \psi(x;\htheta_2)=-u^\top_0 x+1\ge 1-\sqrt{d}\max_{x\in\Omega}\|x\|_\infty|u_0|\ge 0,
\end{align}
so $\psi(\cdot;\htheta)\ge0$ on $\Omega$ for every $\htheta\in\Theta$. Moreover,
\begin{equation}
\int_{\Theta}\psi(x;\htheta)\nu(\td \htheta) = \int_{\{\htheta_1,\htheta_2\}}(u^\top x+v)\nu(\td \htheta)= \tfrac12(u^\top_0 x+1)+\tfrac12(-u^\top_0 x+1)=1,~ \forall x\in\Omega.
\end{equation}
For $\psi(x;\htheta)=u^\top x+v$, $M_{\psi}\leq \sqrt{d}\max_{x\in\Omega}\|x\|_\infty|u|+|v|$ and we have
\begin{multline}
\int_{\Theta} M_{\psi}(\htheta)^2\nu(\td \htheta) =\int_{\{\htheta_1,\htheta_2\}} M_{\psi}(\htheta)^2\nu(\td \htheta) \\\leq \int_{\{\htheta_1,\htheta_2\}} (\sqrt{d}\max_{x\in\Omega}\|x\|_\infty|u|+|v|)^2\nu(\td \htheta) 
=  2(\sqrt{d}\max_{x\in\Omega}\|x\|_\infty|u_0|+1)^2\leq 4.
\end{multline}

As discussed above, the parametrized function $\psi(x;\hat\theta)$ of QTNNs is constructed so that it satisfies the assumptions of Theorem \ref{thm3}. Therefore, for QTNNs, Theorem \ref{thm3} holds. We emphasize that the hypothesis space of parameters $\Theta$ does not need to be finite or consist of only a few specially chosen points. In the examples above, we just use a discrete measure with two atoms merely for clarity. In fact, much richer choices are possible: one may take any symmetric probability measure $\nu$ on a parameter set $\Theta \subset \mbR^d \times \mbR$ as long as its support is contained in a region where $\psi$ satisfies the assumptions of Theorem \ref{thm3}. With such a choice, $\Theta$ can contain many parameters. 

Now we present the approximation property of QTNNs. Following \eqref{08}, the class of two-layer ReLU QTNNs of width $M$ without the output bias is given by
\begin{equation}\label{21}
\mathcal{F}_{M,Q} = 
\Big\{ \phi:\mbR^{d}\rightarrow\mbR~|~\phi(x) = \frac{1}{M} \sum_{m=1}^{M} a_{m} \sigma((w_{m}^\top x + b_{m} )( u_{m}^\top x + v_{m} )), \forall ( a_{m}, w_{m}, b_{m}, u_{m}, v_{m} ) \in \mathcal{D}\Big\}, 
\end{equation}
where $\mathcal{D}=\mbR^{d} \times \mbR \times \Theta$. Noting that
\begin{multline}
M_\psi(\htheta)=\max(\sup_{x\in \Omega}|u^\top x + v|, \sup_{x\in \Omega}|\nabla(u^\top x + v)|) \leq \max(\sup_{x\in\Omega}|u^\top x + v|, |u|)\\
\leq  \max(\sup_{x\in\Omega}(|u| |x| + |v|), |u|) \leq \max(R_\Omega|u| + |v|, |u|) \leq  R_\Omega(|u| + |v|), 
\end{multline}
where $R_\Omega = \max(\sup_{x\in\Omega}|x|,1)$, we can estimate the corresponding Barron-type constant by
\begin{multline}\label{66}
\mathcal{C}_\psi(f) = \left(\inf_{\pi\in\Pi_{f}} \int_{\mathcal{D}} a^2M_{\psi}(\htheta)^2(|w| + |b|)^2\pi(\td a, \td w, \td b, \td \htheta) \right)^{1/2}\\
\leq R_\Omega\left(\inf_{\pi \in \Pi_{f}} \int_{\mathcal{D}}  a^{2}(|w| + |b|)^{2}(|u| + |v|)^{2}  \pi(\td a, \td w, \td b, \td u, \td v)\right)^{1/2} \\
= R_\Omega\left(\inf_{\pi \in \Pi_{f}} \mathbb{E}_{\pi}\left[a^2\left(|w| + |b|\right)^2\left(|u| + |v|\right) ^2\right]\right)^{1/2}\triangleq R_\Omega\mathcal{C}_Q(f),
\end{multline}
where $\Pi_{f}$ is the set of all probability measures $\pi$ such that $$\int_{\mathcal{D}_Q} a\sigma((w^\top x + b)(u^\top x + v)) \pi(\td a, \td w, \td b, \td u, \td v)  = f(x)$$ almost everywhere. Then, by combining Theorem \ref{thm1} and Theorem \ref{thm3}, we obtain the following approximation property for QTNNs.
\begin{corollary}\label{coro1}
Let $\Omega \subset \mbR^{d}$ be a bounded set. Then for any $f\in H^s(\Omega)$ with $s = \lfloor d/2 \rfloor + 3$, there exists some $\phi \in \mathcal{F}_{M,Q}$ such that 
\begin{equation}
\|\phi - f\|_{H^{1}(\Omega)} \leq C_\Omega \mathcal{C}_Q(f) M^{-\frac{1}{2}}, 
\end{equation}
where $C_{\Omega} := \sqrt{10}|\Omega|^{\frac{1}{2}}\max^2(\sup_{x\in \Omega}|x|, 1)$.
\end{corollary}

Corollary \ref{coro1} can be proved using a similar argument to that of Theorem \ref{thm1}, \eqref{66}, and the fact that $\psi(x;\htheta)=u^\top x+v$ satisfies the assumptions of Theorem \ref{thm3}. 

\subsubsection{Application to CTNNs}
The above methods can be adapted for CTNNs. In this section, we turn to the CTNN architecture. Recall that the parametrized function $\psi$ of CTNN is $\psi(x;\theta)=x^\top Qx + u^\top x + v$. We begin with the simplest nontrivial situation in which the parameter hypothesis space $\Theta$ has two elements. Since $\Omega$ is bounded, we denote $R=\max_{x\in\Omega}\|x\|_\infty <\infty$, and can take a matrix $Q_0\in\mathbb{R}^{d\times d}$ and $u_0\in\mathbb{R}^d$ such that 
\begin{equation}
    dR^2\|Q_0\|_2+\sqrt{d}R|u_0|\le 1,
\end{equation}
where $\|\cdot\|_2$ is denoted as the spectral norm of the matrix, and set $\htheta_1=(Q_0,u_0,1)\in\Theta$, $\htheta_2=(-Q_0,-u_0,1)\in\Theta$. We define the probability measure $\nu:=\tfrac12\delta_{\htheta_1}+\tfrac12\delta_{\htheta_2}$, where $\delta_{\htheta}$ is the Dirac mass at $\htheta$. We take $\mathcal{D} = \mbR^{d} \times \mbR \times \Theta $ where $\Theta = \{\htheta_1,\htheta_2\}$.
For $x\in\Omega$ we have
\begin{equation}
|x^\top Q_0x|\le dR^2\|Q_0\|_2,\quad |u^\top_0 x|\le \sqrt{d}R|u_0|.
\end{equation}
Then for all $x\in\Omega$, we obtain
\begin{align}
    &\psi(x;\htheta_1) =x^\top Q_0x + u^\top_0 x + 1 \ge 1-dR^2\|Q_0\|_2-\sqrt{d}R|u_0| \ge 0,\\
    &\psi(x;\htheta_2)=-x^\top Q_0x - u^\top_0 x + 1 \ge 1-dR^2\|Q_0\|_2 -\sqrt{d}R|u_0| \ge 0,
\end{align}
so $\psi(\cdot;\htheta)\ge0$ on $\Omega$ for every $\htheta\in\Theta$. Furthermore,
\begin{multline}
\int_{\Theta}\psi(x;\htheta)\nu(\td \htheta) = \int_{\{\htheta_1,\htheta_2\}}(x^\top Qx + u^\top x + v)\nu(\td \htheta)\\
= \tfrac12\big(x^\top Q_0x + u^\top x + 1\big) +\tfrac12\big(x^\top(-Q_0)x + (-u_0)^\top x + 1\big)
=1,~ \forall x\in\Omega.
\end{multline}
For $\psi=x^\top Qx + u^\top x + v$, $M_{\psi}\leq \max(dR^2\|Q\|_2 + \sqrt{d}R|u|+|v|, 2\sqrt{d}R\|Q\|_2 + |u|)$ and we have
\begin{multline}
\int_{\Theta} M_{\psi}(\htheta)^2 \nu(d\htheta) = \int_{\{\htheta_1,\htheta_2\}} M_{\psi}(\htheta)^2 \nu(d\htheta)  \leq \\
\int_{\{\htheta_1,\htheta_2\}} \max(dR^2\|Q\|_2 + \sqrt{d}R|u|+|v|, 2\sqrt{d}R\|Q\|_2 + |u|)\nu(d\htheta) \\
= 2\max(dR^2\|Q_0\|_2 + \sqrt{d}R|u_0|+1, 2\sqrt{d}R\|Q_0\|_2 + |u_0|)<\infty.
\end{multline}

Therefore, the cubic parametrized function $\psi(x;\hat\theta)$ is chosen so that it also satisfies the assumptions of Theorem \ref{thm3}, and hence Theorem \ref{thm3} applies to CTNNs as well. We stress that this does not impose any restriction that the parameter set $\Theta$ must be finite or contain only a few hand-picked points. In the examples above, we deliberately used a discrete measure with two atoms only to keep the presentation simple. More generally, one may work with any symmetric probability measure $\nu$ supported on a set $\Theta \subset \mathbb{R}^{d\times d}\times \mathbb{R}^d \times \mathbb{R}$, as long as the support of $\nu$ lies in a region where $\psi$ satisfies the assumptions of Theorem \ref{thm3}. In this way, $\Theta$ can encode a genuinely large and flexible family of parameters, rather than being restricted to a few specific examples.

Specifically, from \eqref{07}, the class of two-layer ReLU CTNNs of width $M$ is given by
\begin{equation}
\mathcal{F}_{M,C} = 
\biggl\{ \phi:\mbR^{d}\rightarrow\mbR~|~\phi(x) = \frac{1}{M} \sum_{i=1}^{M} a_{i} \sigma\bigl[( w_{i}^\top x + b_{i} )(x^\top Q_i x + u_{i}^\top  x + v_{i} )\bigr],  \forall ( a_{i}, w_{i}, b_{i}, Q_i, u_{i}, v_{i} ) \in \mathcal{D} \biggr\}.
\end{equation}
where $\mathcal{D}=\mbR^{d} \times \mbR \times \Theta$. Recall that
\begin{multline}
M_\psi(\htheta)=\max\left(\sup_{x\in \Omega}|x^\top Qx + u^\top x + v|, \sup_{x\in \Omega}|\nabla(x^\top Qx + u^\top x + v)|\right) \\ \leq \max\left(\sup_{x\in\Omega}|x^\top Qx + u^\top x + v|, |Qx + u|\right)
\leq  \max\left(\sup_{x\in\Omega}(\|Q\|_2 |x|^2 + |u||x| + |v|), \sup_{x\in\Omega}(\|Q\|_2|x| + |u|)\right) \\
\leq \max\left(R^2_\Omega\|Q\|_2 + R_\Omega|u| + |v|, R_\Omega\|Q\|_2 + |u|\right) \leq  R^2_\Omega(\|Q\|_2 + |u| + |v|), 
\end{multline}
where $R_\Omega = \max(\sup_{x\in\Omega}|x|,1)$. We can estimate the corresponding Barron-type constant by
\begin{multline}\label{67}
\mathcal{C}_\psi(f) = \left(\inf_{\pi\in\Pi_{f}} \int_{\mathcal{D}} a^2M_{\psi}(\htheta)^2(|w| + |b|)^2\pi(\td a, \td w, \td b, \td\htheta) \right)^{1/2}\\
\leq R^2_\Omega\left(\inf_{\pi \in \Pi_{f}} \int_{\mathcal{D}}  a^2(\|Q\|_2 + |u| + |v|)^2(|w| + |b|)^2  \pi(\td a, \td w, \td b,\td Q, \td u, \td v)\right)^{1/2} \\
= R^2_\Omega\left(\inf_{\pi \in \Pi_{f}} \mathbb{E}_{\pi}\left[a^2\left(\|Q\|_2 + |u| + |v|\right)^2\left(|w| + |b|\right) ^2\right]\right)^{1/2}\triangleq R^2_\Omega\mathcal{C}_C(f),
\end{multline}
where $\Pi_{f}$ is the set of all probability measures $\pi$ such that 
\begin{equation}
    \int_{\mathcal{D}} a\,\sigma\left((w^\top x + b)(x^\top Qx+u^\top x+v)\right) \pi(\td a, \td w, \td b,\td Q, \td u, \td v)  = f(x),
\end{equation}
almost everywhere. Then, by combining Theorem \ref{thm1}, Theorem \ref{thm3} and \eqref{67}, we have the following approximation property for CTNNs.
\begin{corollary}\label{coro4}
Let $\Omega \subset \mbR^{d}$ be a bounded set. Then for any $f\in H^s(\Omega)$ with $s = \lfloor d/2 \rfloor + 3$, there exists some $\phi \in \mathcal{F}_{M,C}$ such that 
\begin{equation}
\|\phi - f\|_{H^{1}(\Omega)} \leq C_\Omega \mathcal{C}_C(f) M^{-\frac{1}{2}},
\end{equation}
where $C_{\Omega} := \sqrt{20}|\Omega|^{\frac{1}{2}}\max^3(\sup_{x\in \Omega}|x|, 1)$.
\end{corollary}

\subsection{Error bounds of deep architectures} 
In this section, we derive error bounds of GTNNs with deep architectures. Specifically, we investigate the approximation property of deep GTNNs defined by \eqref{09} with ReLU activation for the smooth function space $C^s([0,1]^d)$. Recall the norm $\|\cdot\|_{C^{s}\left([0,1]^{d}\right)}$ is defined by
\begin{equation}
\|f\|_{C^{s}\left([0,1]^{d}\right)} := 
\max\left\{ 
    \|\partial^{\alpha} f\|_{L^{\infty}\left([0,1]^{d}\right)} : 
    \|\alpha\|_{1} \leq s, \alpha \in \mathbb{N}^{d}
\right\}
,\quad\forall f \in C^s([0,1]^d).
\end{equation}  
Then, the result is given by
\begin{thm} \label{approthm2}
Given a smooth function $f \in C^s([0,1]^d)$ with $s\in \mathbb{N}^+$. Suppose $\psi_\ell(\cdot;\htheta_{\ell,0})$ is a non-zero constant for some $\htheta_{\ell,0}\in\Theta$ in \eqref{09}. Then for any $N,K\in \mathbb{N}^+$, there exists a deep ReLU GTNN $\phi$ \eqref{09} with width $C_1(N+2)\operatorname{log}_2(8N)$ and depth $C_2(K+2)\operatorname{log}_2(4K) + 2d$ such that
\begin{equation}\label{25}
\|\phi - f\|_{L^{\infty}\left([0,1]^{d}\right)} \leq C_{3}\|f\|_{C^{s}\left([0,1]^{d}\right)} N^{-2s/d}K^{-2s/d},
\end{equation}
where $C_1 = 17s^{d+1}3^dd , C_2 = 18s^2$ and $C_3 = 85(s + 1)^d8^s$. 
\end{thm}
\begin{proof}
By Theorem 1.1 in \citep{Lu2021}, there exists a deep fully connected ReLU neural network $\phi$ with width $C_1(N+2)\log_2(8N)$ and depth $C_2(K+2)\log_2(4K) + 2d$, where $C_1 = 17s^{d+1}3^dd , C_2 = 18s^2$ and $C_3 = 85(s + 1)^d8^s$, such that \eqref{25} holds.

Now we consider the deep GTNNs. Taking $\htheta_\ell=\htheta_{\ell,0}$ in \eqref{09}, then $\psi_\ell(x;\htheta_\ell)$ becomes a constant vector-valued function, and the corresponding GTNN $\phi$ is simplified to a fully connected ReLU neural network, which demonstrates that the proof is completed.
\end{proof}

Specifically, Theorem \ref{approthm2} holds for deep QTNNs, CTNNs, and TTNNs. It suffices to make $\psi$ simplified to a constant. For the deep QTNNs, we can take $U_\ell$ as the zero matrix and $v_\ell=[1~\dots~1]^\top$ for all $\ell$ in \eqref{23}; for the deep CTNNs, we additionally take $Q_\ell^i$ ($i=1,\dots, M_\ell$) as zero matrices for all $\ell$ in \eqref{24}; for the deep TTNNs, we take $U^n_\ell$ as the zero matrix and $v^n_\ell=[1~\dots~1]^\top$ for all $\ell$ and $n=1,\dots N$ in \eqref{52}.  

Theorem \ref{approthm2} characterizes the approximation error in terms of the total number of general neurons (with an arbitrary distribution in width and depth) and the smoothness of the target function to be approximated. If $f$ is a $C^s$ smooth function in $[0,1]^d$ with $s\in \mathbb{N}^+$, then there exists a deep ReLU GTNN $\phi$ with width $O(N\ln N)$ and depth $O(K\ln K)$ such that it can approximate $f$ with an approximation error $O(\|f\|_{C^{s}\left([0,1]^{d}\right)} N^{-2s/d}K^{-2s/d})$. This estimate is non-asymptotic in the sense that it is valid for arbitrary width and depth specified by $N\in\mathbb{N}^+$ and $K\in \mathbb{N}^+$, respectively \citep{Lu2021}.

In summary, we have investigated the approximation properties of deep GTNNs. The error bound in Theorem~\ref{approthm2} is explicitly controlled by the width and depth of the network. Beyond the general framework, we have also analyzed three specific GTNN architectures: QTNNs, CTNNs, and TTNNs. Further approximation properties for deep GTNNs can be derived using existing approximation theory of deep neural networks.

\section{Numerical Examples}\label{sec_numerical_examples}
The approximation properties of the proposed QTNNs, CTNNs, and TTNNs are tested by several numerical examples in this section. We test the two-layer ($L=2$) QTNNs, CTNNs and TTNNs defined by \eqref{08}, \eqref{07} and \eqref{28}, respectively. For simplicity, it should be noted that in TTNNs, we only take the $\sin$ term for $\psi$, but this is not limited to this form. And we also test the three-layer $(L=3)$ QTNNs, CTNNs, and TTNNs without bias in the output layer; from \eqref{09}, they are formulated by
\begin{equation}\label{15}
\phi(x) = \sum_{i=1}^M \alpha_{i} \sigma( \widetilde{\phi}_{i}(x) + \beta_{i}),
\end{equation}
where $M \in \mathbb{N}^{+}$ is the width (we use the same width for both hidden layers); $\alpha_{i}\in\mbR$ is the weight of the output layer; $\beta_{i}\in\mbR$ is the bias of the second-layer neurons; $\widetilde{\phi}_{i}$ is chosen as \eqref{08} for the QTNN or chosen as \eqref{07} for the CTNN or chosen as \eqref{28} for the TTNN. Furthermore, we test the two and three-layer FNNs \eqref{04} for comparison. 
For different cases, let the target function be $f(x)$, and we extract $N$ data pairs $(x_n,f(x_n))$ from the target function as the training points. Precisely, we solve the following least squares regression
\begin{equation}\label{14}
\min_{\theta}\frac{1}{N}\sum_{n=1}^N|\phi(x_n,\theta)-f(x_n)|^2,
\end{equation}
where $\phi(x_n,\theta)$ is the learner.

\subsection{Settings}
\begin{itemize}
\item \textbf{Environment.} The methods are implemented mainly in the Python environment. The PyTorch library, combined with the CUDA toolkit, is utilized for implementing neural networks and enabling GPU-based parallel computing.
\item \textbf{Optimizer and hyper-parameters.} The optimization in the experiments is solved by the Adam optimizer \citep{Kingma2015}. The algorithm is implemented for a total of 20000 iterations with learning rates decaying linearly from $10^{-2}$ to $10^{-3}$.
\item \textbf{Network setting.} The activation function $\sigma$ is set as ReLU for regression problems \eqref{14}. We test the networks with depths $L=2,3$. For a fair comparison, for every $L$, we adjust the widths $M$ so that the QTNN/CTNN/TTNN and FNN architectures have approximately the same number of trainable parameters. The parameters of FNNs are initialized by
\begin{equation}
a,W_{l},b_{l}\sim U(-\frac{1}{\sqrt{M}},\frac{1}{\sqrt{M}}),\quad l=1,\ldots, L.
\end{equation}
Similar random initialization under a uniform distribution is used for QTNNs/CTNNs/TTNNs.
\item \textbf{Testing set and error evaluation.} We generate a testing set $\mathcal{X}$ consisting of $10000$ uniform random points in the domain. For the regression problems \eqref{14}, we compute the following relative $\ell^{2}$ error over $\mathcal{X}$:
\begin{equation}\label{27}
e_{\ell^{2}(\mathcal{X})}:=\left(\sum_{x\in\mathcal{X}}|\phi(x)-f(x)|^{2}/\sum_{x\in\mathcal{X}}|f(x)|^{2}\right)^{\frac{1}{2}},
\end{equation}
where $\phi$ is the trained neural network and $f$ is the target function.
\item \textbf{Repeated experiments.} We repeat each experiment with 10 different random seeds, and present the best result among the 10 random tests. The repeated experiments reduce the effect of randomness in initialization and stochastic training, making the results fair and convincing.
\end{itemize}

\subsection{Case 1: Learning oscillatory functions}
To test the performance of GTNNs, we choose the following two oscillatory functions as target functions:
\begin{align}\label{17}
&\text{Case 1.1: 1-D oscillatory function}~f_1(x) = \sin(30\pi x),~ x\in[-1,1],  \\  
&\text{Case 1.2: 5-D oscillatory function}~f_2(x) = \sin(8\pi (|x|^2-1)),~ x \in \mbR^{5},~|x|\leq 1. 
\end{align}

For Case 1.1 and Case 1.2, we select 20,000 training feature points $\{x_n\}$ that are evenly distributed within their domain. We fit the dataset $\{x_n, f_k(x_n)\}$ ($k=1,2$) using CTNNs, QTNNs, TTNNs and FNNs through the least squares regression \eqref{14}.

\begin{figure}[htbp]
\centering
\includegraphics[scale=0.45]{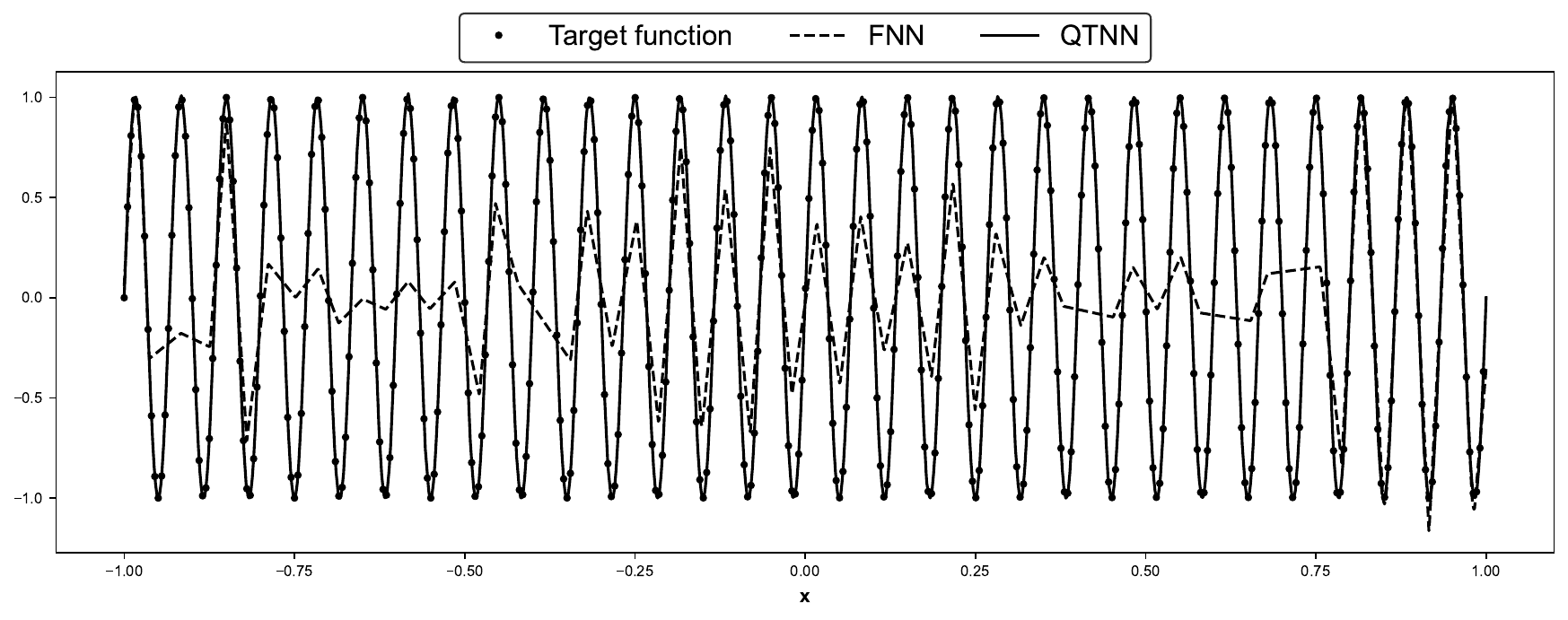}\hspace{10pt}
\caption{\em The target function, the trained FNN and QTNN in Case 1.1.}
\label{fig:case2}
\end{figure}

The function $f_1(x)$ is oscillatory with $30$ waves in the domain, which is usually hard to learn using standard neural networks. It is observed in Table \ref{tab:case2.1} that FNNs can only approximate $f_1(x)$ with error up to $O(10^{-1})$ with $L=3$, implying FNNs cannot effectively learn such high-frequency functions. Comparatively, QTNNs can obtain lower errors $O(10^{-3})$ with $L=3$ and CTNNs/TTNNs can obtain lower errors $O(10^{-2})$ with $L=3$ (see Table \ref{tab:case2.1}). Moreover, we plot the target function $f_1(x)$, the trained FNN learner, and the trained QTNN learner in Figure \ref{fig:case2}, and it is clear that QTNNs are more capable of learning the high-frequency waves than FNNs. 

\begin{table}[htbp]
\footnotesize
\centering
\caption{Relative $\ell^2$ errors of CTNN, QTNN, TTNN and FNN in Case 1.1.}
\label{tab:case2.1}
\setlength{\tabcolsep}{6pt}
\renewcommand{\arraystretch}{1.15}
\begin{tabular}{c l cccc}
\toprule
$L$ & Parameters & CTNN & QTNN & TTNN & FNN \\
\midrule
\multirow{4}{*}{$2$}
 & $\approx100$ & 9.5670e-01 & 9.3983e-01 & 9.6293e-01 & 9.9184e-01 \\
 & $\approx300$ & 6.7881e-01 & 7.8244e-01 & 8.2942e-01 & 9.8804e-01 \\
 & $\approx400$ & 5.6866e-01 & 6.8654e-01 & 7.7628e-01 & 9.8862e-01 \\
 & $\approx500$ & \textbf{4.8583e-01} & 6.5233e-01 & 6.6346e-01 & 9.8242e-01 \\
\midrule
\multirow{4}{*}{$3$}
 & $\approx2000$  & 2.3849e-01 & 2.2851e-01 & 1.3483e-01 & 9.3881e-01 \\
 & $\approx18000$ & 1.8402e-02 & 1.7308e-02 & 2.3252e-02 & 8.3004e-01 \\
 & $\approx32000$ & 1.2788e-02 & 9.3258e-03 & 1.4913e-02 & 7.7385e-01 \\
 & $\approx50000$ & 8.9352e-03 & \textbf{6.5315e-03} & 1.2822e-02 & 7.5395e-01 \\
\bottomrule
\end{tabular}
\end{table}

\begin{figure}[htbp]
\centering
\includegraphics[scale=0.45]{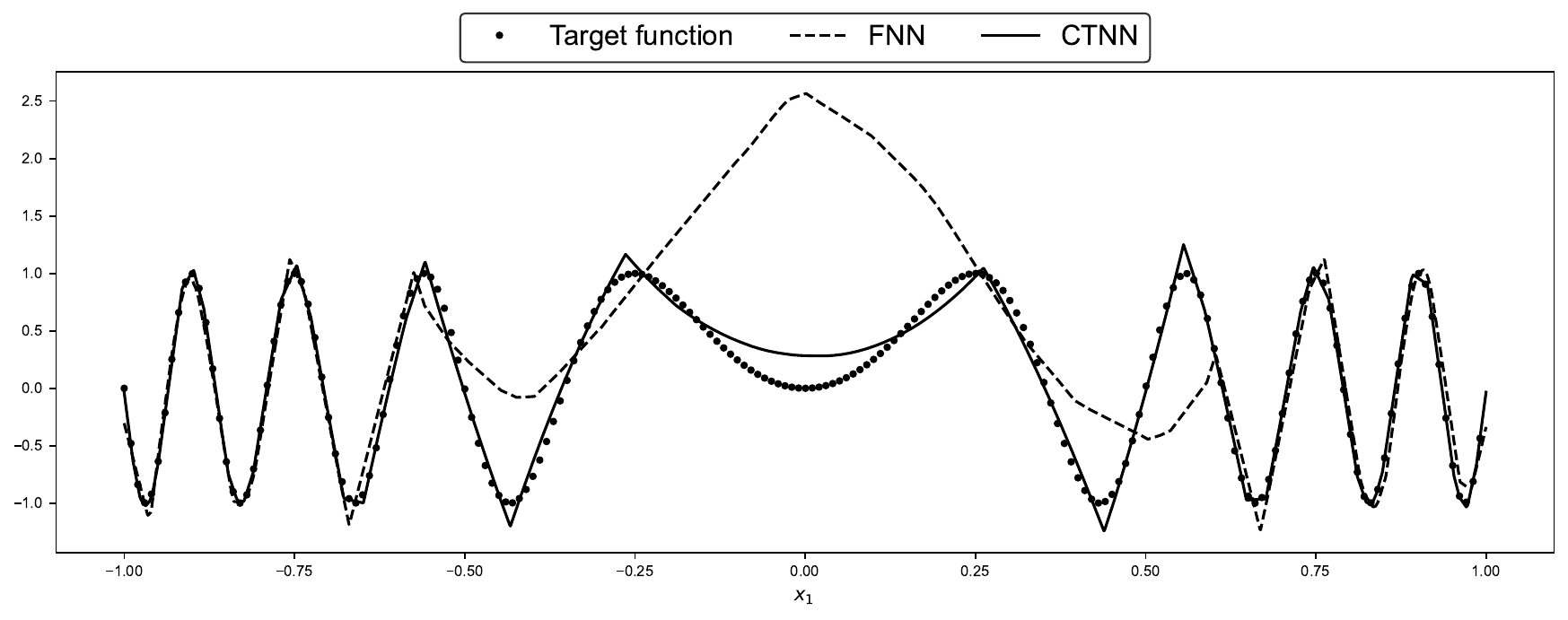}
\caption{\em Graphs of the target function and the trained FNN and CTNN at $(x_1,0,0,0,0)$ in Case 1.2.}
\label{fig:case2.3_1}
\end{figure}

\begin{figure}[htbp]
\centering
\subfloat[target function]{
\includegraphics[scale=0.3, trim=100 200 100 230,clip]{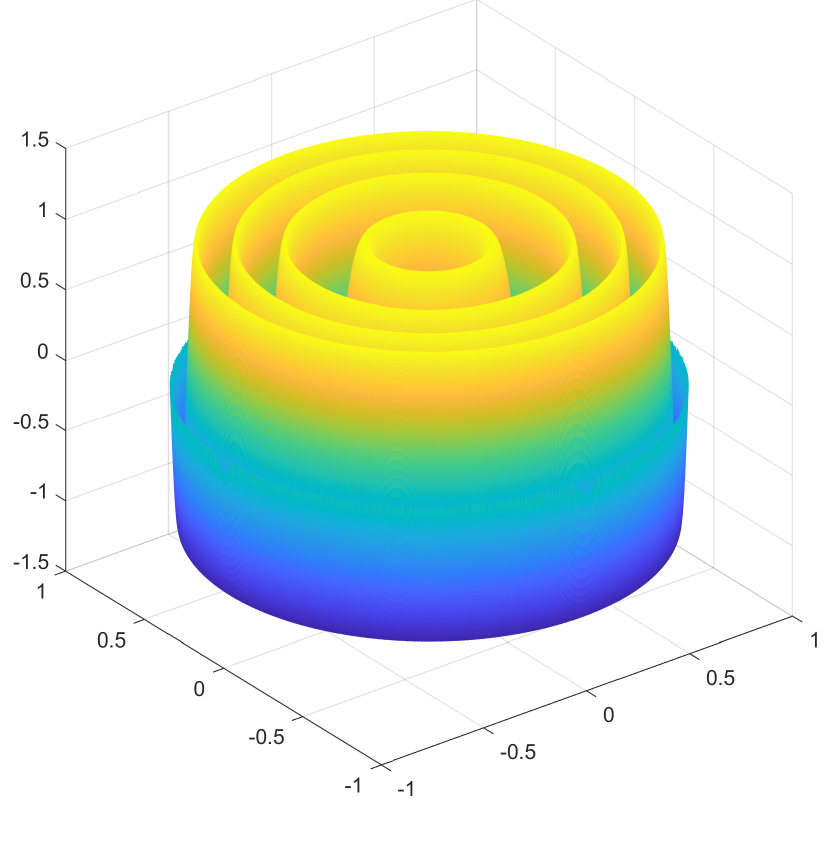}}\hspace{10pt}
\subfloat[CTNN]{
\includegraphics[scale=0.3, trim=100 200 100 230,clip]{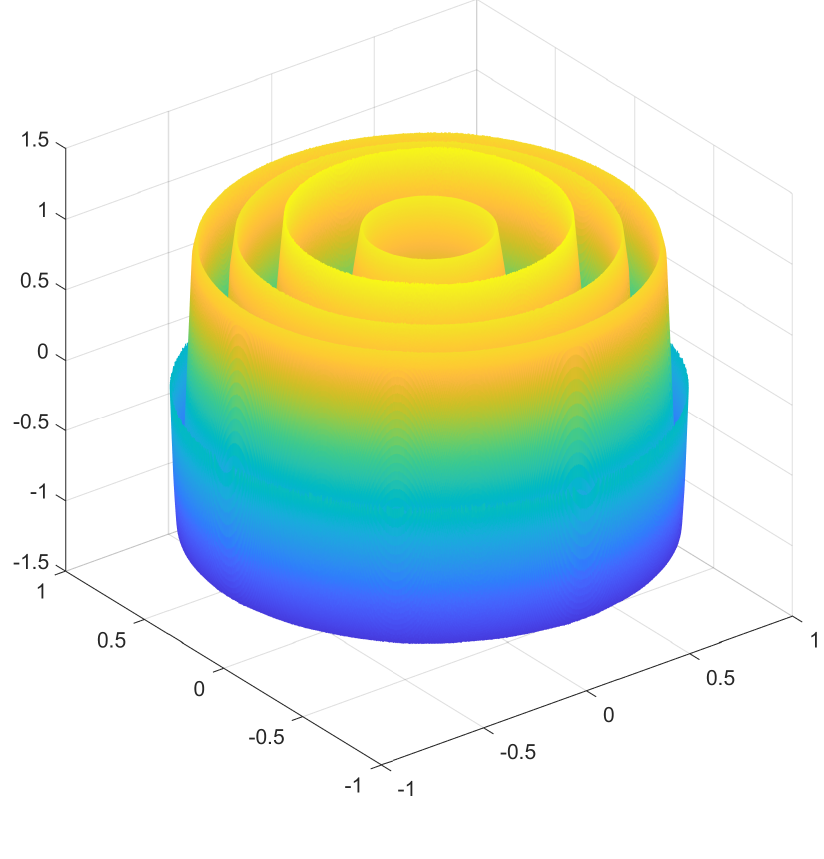}}\hspace{10pt}
\subfloat[FNN]{
\includegraphics[scale=0.3, trim=100 200 100 230,clip]{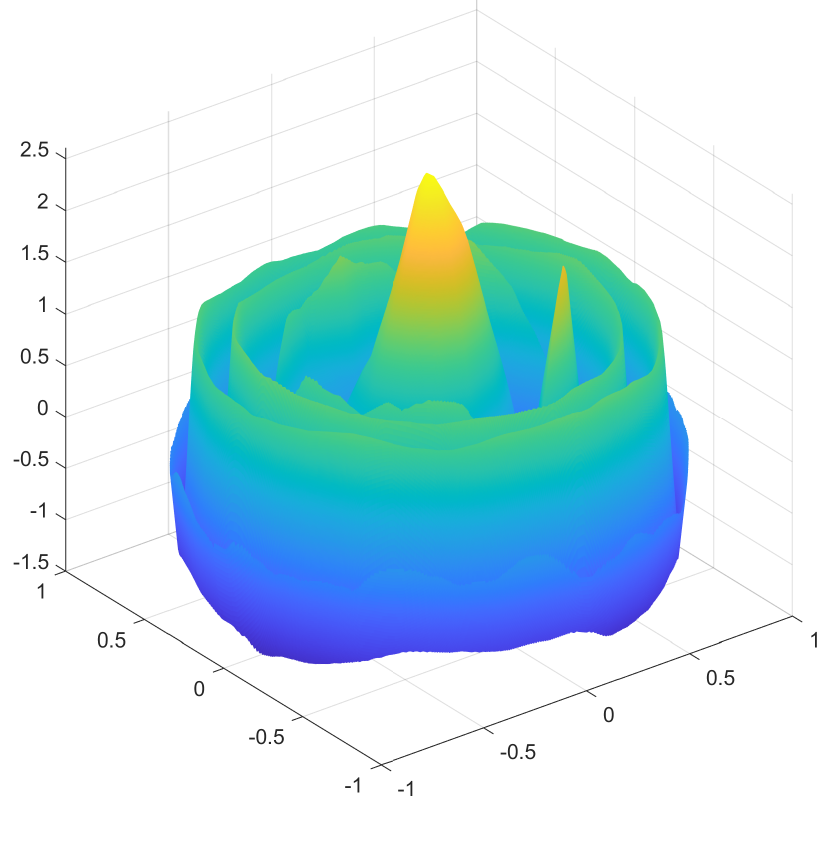}}
\caption{\em 2-D profiles of the target function and the trained FNN and CTNN at $(x_1,x_2,0,0,0)$ in Case 1.2.}
\label{fig:case2.3_2}
\end{figure}

The function $f_2(x)$ is of 5 dimensions, which is typically difficult to learn using linear structures such as polynomials or finite elements. Here, we test the learning capability of GTNNs on these slightly high-dimensional functions. Moreover, we show the 1-D and 2-D profiles of $f_2(x)$ in Figure \ref{fig:case2.3_1} and \ref{fig:case2.3_2}, respectively. Note that $f_2$ is oscillatory with eight waves visible in its profile. It is observed that CTNN can roughly capture all waves, but FNN completely fails in learning the wave near the origin. The errors are listed in Table \ref{tab:case2.3}, and it is seen that CTNN achieves $O(10^{-2})$ errors, but FNN/QTNN/TTNN can only obtain $O(10^{-1})$ errors. This indicates that CTNNs can learn this highly oscillatory, high-dimensional function much more effectively than QTNNs, TTNNs and FNNs under the same parameter budget. 

\begin{table}[htbp]
\footnotesize
\centering
\caption{Relative $\ell^2$ errors of CTNN, QTNN, TTNN and FNN in Case 1.2.}
\label{tab:case2.3}
\setlength{\tabcolsep}{6pt}
\renewcommand{\arraystretch}{1.15}
\begin{tabular}{c l cccc}
\toprule
$L$ & Parameters & CTNN & QTNN & TTNN & FNN \\
\midrule
\multirow{4}{*}{\textit{$2$}}
 & $\approx750$  & 9.4215e-02 & 9.6595e-01 & 9.6514e-01 & 9.6947e-01 \\ 
 & $\approx2300$ & 6.8720e-02 & 9.9279e-01 & 9.9514e-01 & 9.8577e-01 \\ 
 & $\approx3000$ & 6.6818e-02 & 1.0042e+00 & 1.0129e+00 & 9.9113e-01 \\ 
 & $\approx3800$ & \textbf{5.9711e-02} & 1.0225e+00 & 1.0264e+00 & 1.0008e+00 \\ 
\midrule
\multirow{4}{*}{\textit{$3$}}
 & $\approx1000$  & 2.4484e-01 & 8.9331e-01 & 8.8925e-01 & 6.6103e-01 \\ 
 & $\approx4000$  & 5.2547e-02 & 3.0307e-01 & 5.3282e-01 & 3.5751e-01 \\ 
 & $\approx9000$  & 4.9712e-02 & 1.5137e-01 & 3.4199e-01 & 3.0405e-01 \\ 
 & $\approx15000$ & \textbf{4.5762e-02} & 1.2416e-01 & 2.5668e-01 & 2.5744e-01 \\ 
\bottomrule
\end{tabular}
\end{table}

\subsection{Case 2: Learning the flow of the shock-turbulence problem}
We consider the flow of the shock-turbulence problem \citep{Titarev2004} on the interval $[-1, 1]$, which exhibits oscillatory behavior with numerous high-frequency waves. We extract $N=2000$ data pairs $(x_n,y_n)$ from the flow, and use CTNNs, QTNNs, TTNNs and FNNs to fit the data. Precisely, we solve the least squares optimization \eqref{14} by replacing $f(x_n)$ with $y_n$. 

\begin{figure}[htbp]
\centering
\includegraphics[scale=0.45]{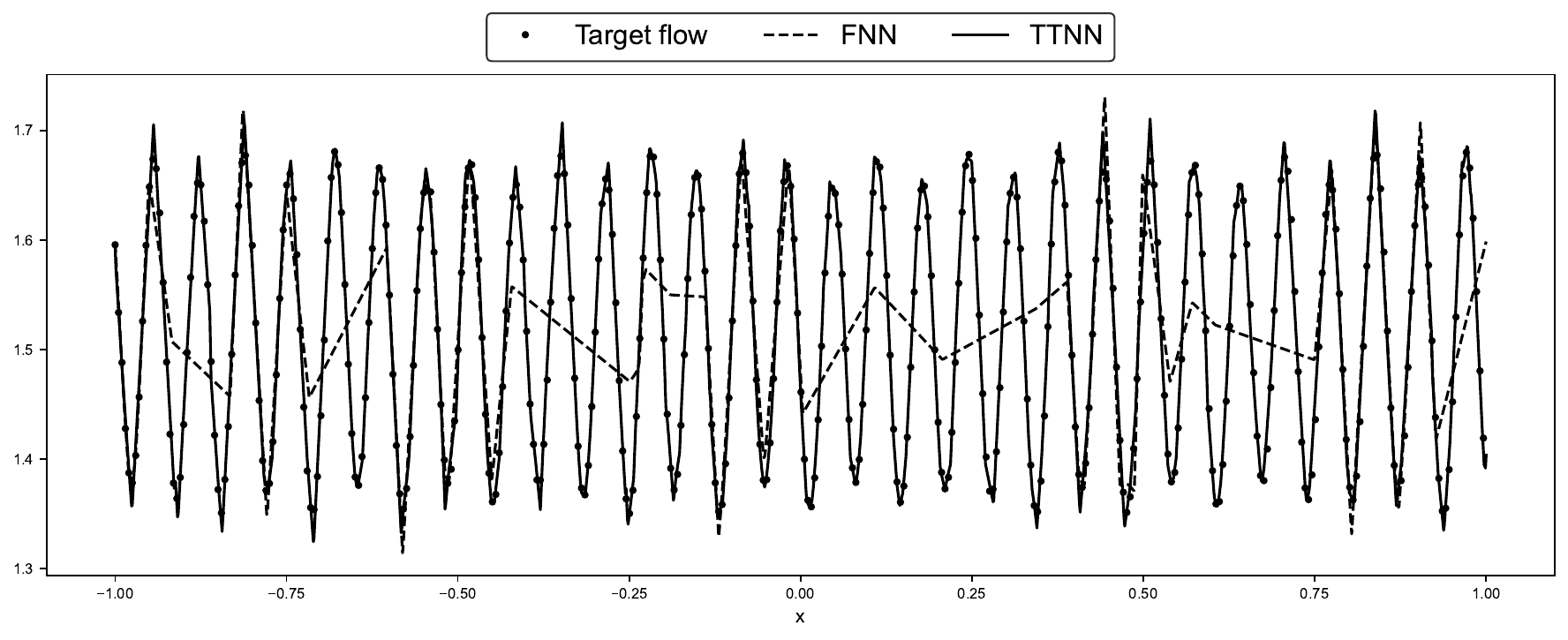}
\caption{\em The target flow, the trained FNN and the trained TTNN in Case 2.}
\label{fig:case1}
\end{figure}

The testing errors for various numbers of parameters and depth $L$ are presented in Table \ref{tab:case1}. It is observed that CTNNs/QTNNs/TTNNs can obtain errors of $O(10^{-3})$, but FNNs have larger errors of $O(10^{-2})$. When $L=3$, the benefit of the GTNN architectures becomes much more pronounced. As the number of parameters grows from $600$ to $21000$, the errors of CTNNs, QTNNs, and TTNNs decrease steadily from $O(10^{-2})$ to $O(10^{-3})$, whereas FNNs have larger errors of $O(10^{-2})$. Moreover, we present the shapes of the target flow, the FNN learner, and the TTNN learner in Figure \ref{fig:case1}. Here, the TTNNs can effectively learn the high-frequency features of the data, whereas FNNs struggle to learn part of the oscillating waves. 
\begin{table}[htbp]
\footnotesize
\centering
\caption{Relative $\ell^2$ errors of CTNN, QTNN, TTNN and FNN in Case 2.}
\label{tab:case1}
\setlength{\tabcolsep}{6pt}
\renewcommand{\arraystretch}{1.15}
\begin{tabular}{c l cccc}
\toprule
$L$ & Parameters & CTNN & QTNN & TTNN & FNN \\
\midrule
\multirow{4}{*}{\textit{$2$}}
 & $\approx 200$  & 6.5994e-02 & 6.6016e-02 & 6.5293e-02 & 6.8479e-02 \\
 & $\approx 400$  & 4.8571e-02 & 5.9492e-02 & 5.6068e-02 & 6.5279e-02 \\
 & $\approx 600$  & 4.8571e-02 & 4.9449e-02 & 4.9914e-02 & 6.3550e-02 \\
 & $\approx 1000$ & \textbf{3.4737e-02} & 4.6596e-02 & 4.5540e-02 & 6.0641e-02 \\
\midrule
\multirow{4}{*}{\textit{$3$}}
 & $\approx 600$   & 4.4816e-02 & 5.9467e-02 & 5.3806e-02 & 6.9865e-02 \\
 & $\approx 2500$  & 2.9431e-02 & 3.5746e-02 & 2.8723e-02 & 6.4602e-02 \\
 & $\approx 10000$ & 6.2446e-03 & 1.0257e-02 & 5.2829e-03 & 5.9516e-02 \\
 & $\approx 21000$ & 4.8819e-03 & 4.9412e-03 & \textbf{4.8122e-03} & 5.3214e-02 \\
\bottomrule
\end{tabular}
\end{table}

\subsection{Case 3: Learning the multi-peak Gaussian functions}
To further test the ability of GTNNs in learning localized features with multiple sharp peaks, we consider the following Gaussian-mixture type target functions:
\begin{align}
&\text{Case 3.1: 1-D multi-Gaussian function}\quad f_3(x)=\sum_{i=1}^{7} a_i e^{-\frac{(x-c_i)^2}{2s_i^2}},~ x\in[-1,1];\\
&\text{Case 3.2: 2-D multi-Gaussian function}\quad f_4(x)=\sum_{i=1}^{10}\alpha_i e^{-\frac{\|x-\mu_i\|_2^2}{2\sigma_i^2}},~ x\in[-1,1]^2.
\end{align}
For Case 3.1, the parameters are fixed as
\begin{align*}
&(c_1,c_2,c_3,c_4,c_5,c_6)=(-0.90,-0.52,-0.16,0.02,0.38,0.82,0.85),\\
&(s_1,s_2,s_3,s_4,s_5,s_6)=(0.030,0.032,0.030,0.024,0.027,0.025,0.6),\\
&(a_1,a_2,a_3,a_4,a_5,a_6)=(1.00,1.10,1.20,0.90,0.88,0.92,0.1).
\end{align*}
For Case 3.2, the amplitudes and standard deviations are sampled as $\alpha_i \sim U(0.85,\,1.35)$,  $\sigma_i \sim U(0.11,\,0.20)$, respectively and the centers satisfy $\mu_i=(c_{i1},c_{i2})$ with $c_{i1},c_{i2}\sim U(-0.92,\,0.92)$. We sample $20000$ input points $\{x_n\}$ evenly across their domain and fit the data pairs $\{(x_n,f_k(x_n))\}$ using CTNNs, QTNNs, TTNNs, and FNNs with the least-squares regression \eqref{14}.

\begin{figure}[htbp]
\centering
\includegraphics[scale=0.6]{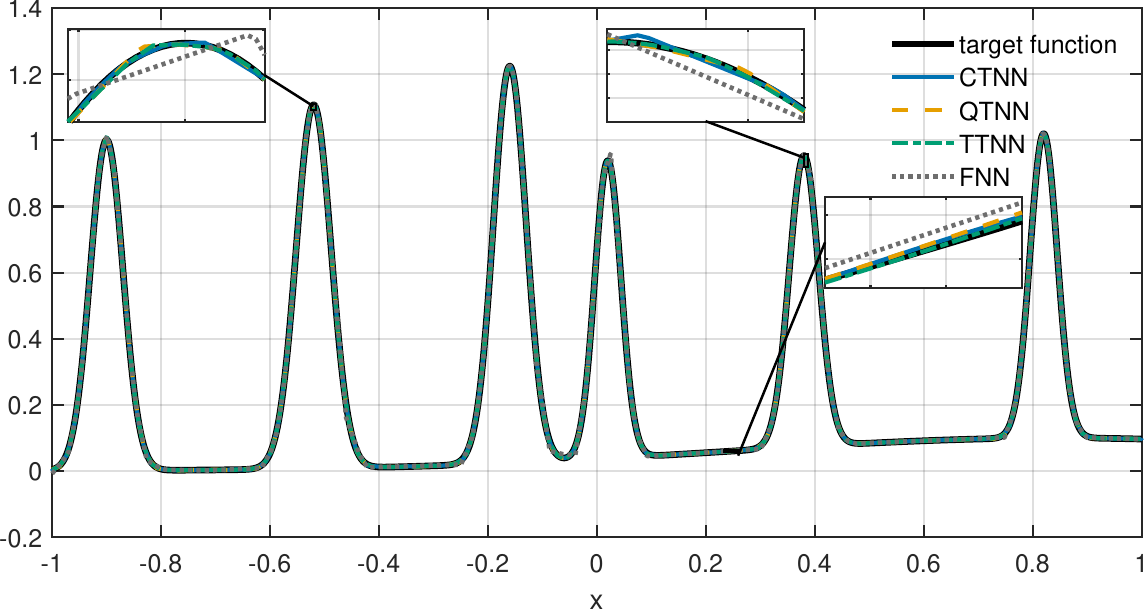}
\caption{\em The target function, the trained FNN, CTNN, QTNN and TTNN in Case 3.1.}
\label{fig:case3.1}
\end{figure}

The 1-D multi-peak Gaussian target $f_3(x)$ is essentially a ``local-feature" function: it consists of several very narrow peaks separated by long, almost flat regions, together with a mild background bump. Figure \ref{fig:case3.1} shows that the trained CTNN, QTNN, and TTNN track the curve of $f_3(x)$ across the whole interval, and the agreement remains tight in the zoomed-in views, where small discrepancies are easiest to spot. By contrast, the FNN exhibits a visible bias in these magnified regions, especially along the shoulders and tails of the peaks, suggesting that it tends to smooth out or slightly shift the localized structures under the same parameter budget.

The quantitative results in Table \ref{tab:case3.1} support these observations. For $L=2$, CTNNs/QTNNs/TTNNs already reach the $O(10^{-2})$ error level with only a few hundred parameters, while the FNN remains at the $O(10^{-1})$ level. When $L=3$, the errors of CTNNs/QTNNs/TTNNs further decrease to the $O(10^{-3})$, whereas the FNN still stays around $O(10^{-2})$. Overall, this case highlights a clear advantage of GTNNs on functions dominated by sharp, localized peaks: they achieve noticeably better accuracy with the same number of parameters as FNN.

\begin{figure}[tbh!]
\centering
\subfloat[target function $f_4(x)$]{
\includegraphics[scale=0.45]{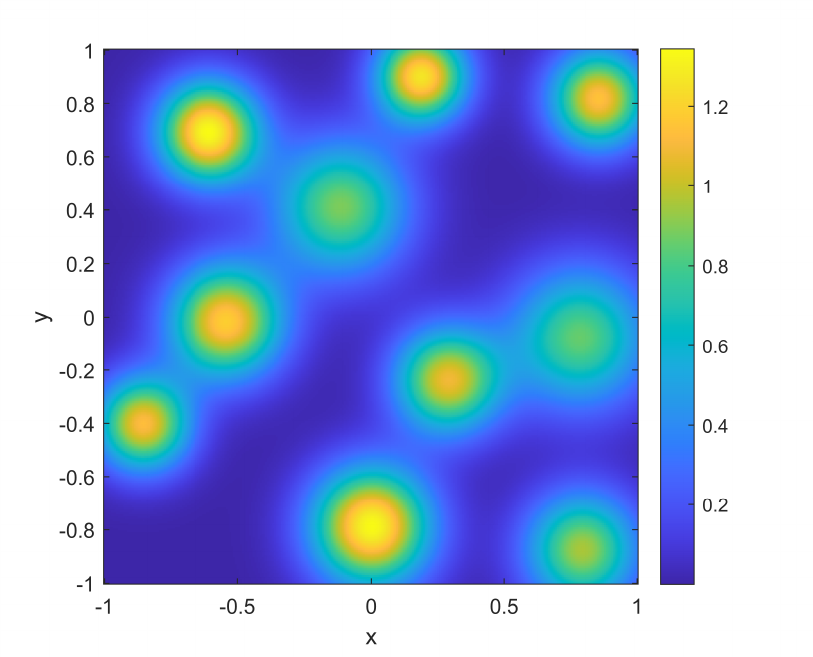}}\hspace{10pt}
\subfloat[$|\phi_{FNN}-f_4(x)|$]{
\includegraphics[scale=0.45]{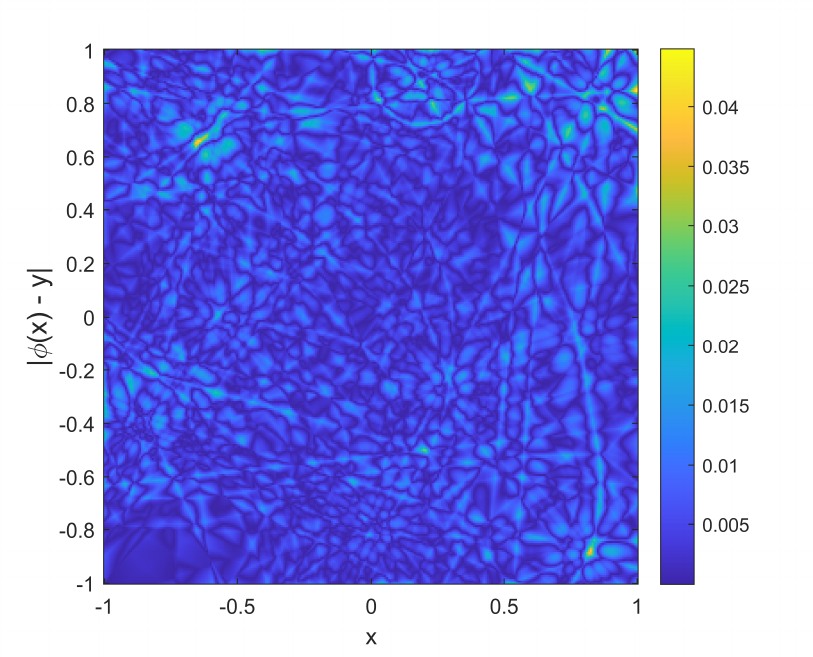}}\hspace{10pt}
\subfloat[$|\phi_{TTNN}-f_4(x)|$]{
\includegraphics[scale=0.45]{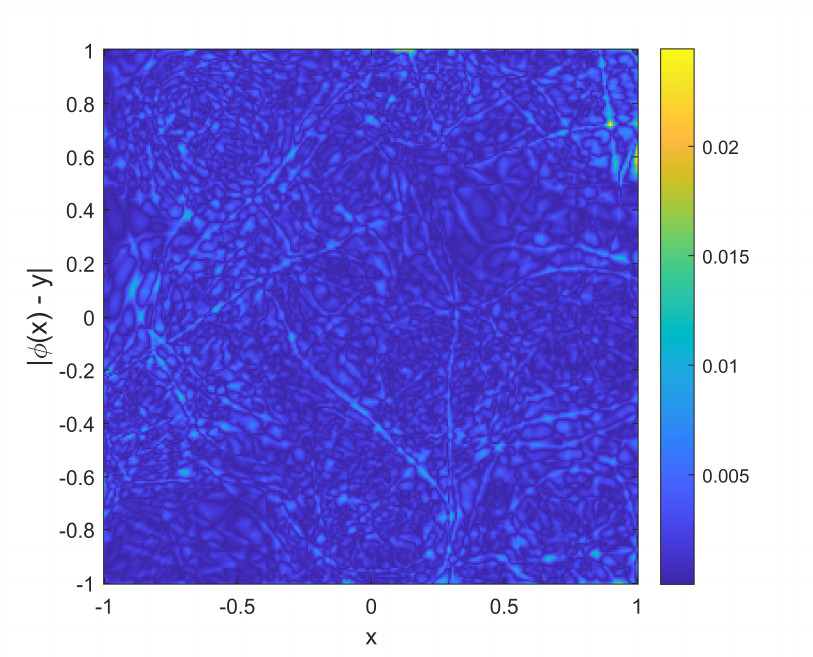}}\hspace{10pt}
\subfloat[$|\phi_{CTNN}-f_4(x)|$]{
\includegraphics[scale=0.45]{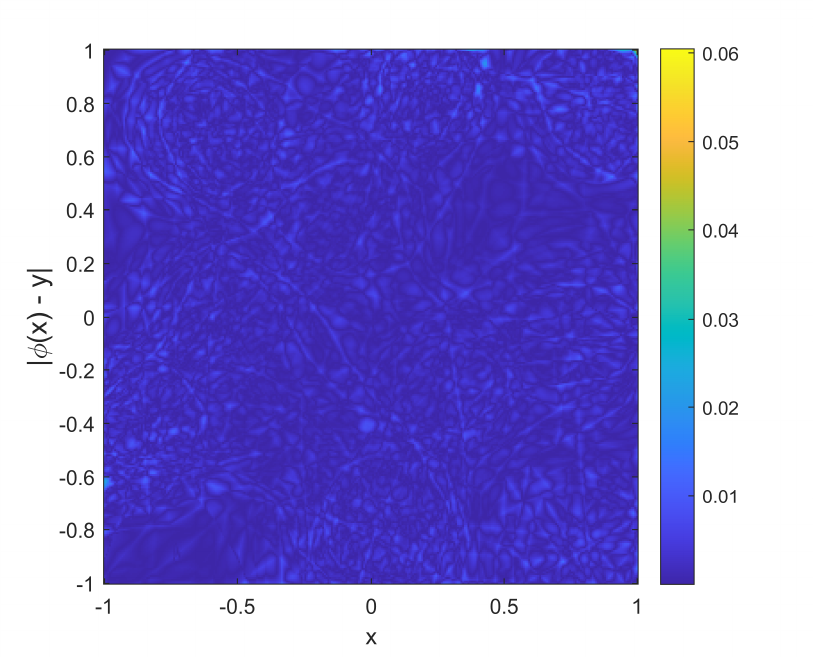}}\hspace{10pt}
\subfloat[$|\phi_{QTNN}-f_4(x)|$]{
\includegraphics[scale=0.45]{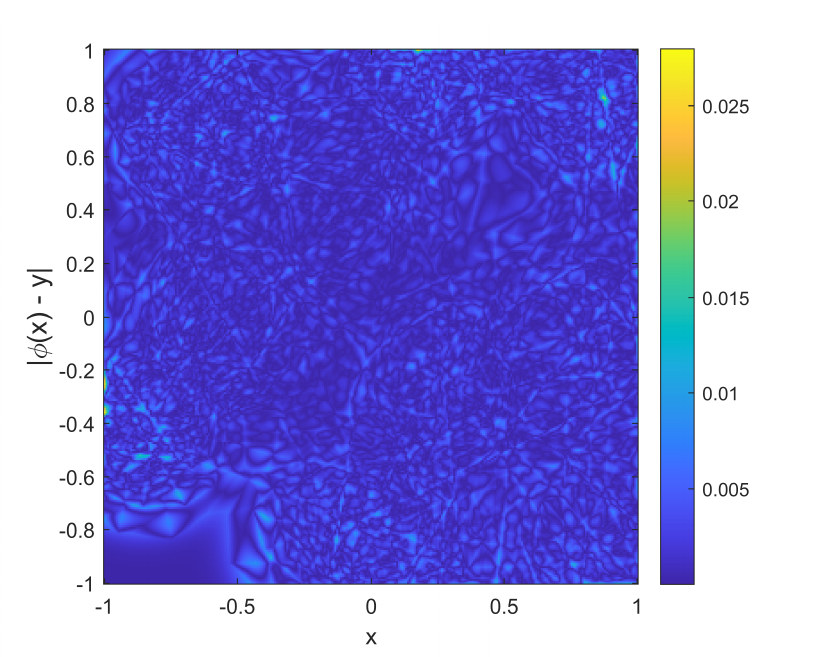}}
\caption{\em 2-D profiles of the target function and absolute errors of the trained FNN, QTNN, CTNN and TTNN in Case 3.2.}
\label{fig:case3.2}
\end{figure}

As shown in Figure \ref{fig:case3.2}, the TTNNs/CTNNs/QTNNs match $f_4(x)$ very closely and preserve the shapes of the peaks. In comparison, the FNN result is slightly more diffusive around peak shoulders, where the function transitions most rapidly. The FNN errors exhibit more pronounced patches near steep transition regions, whereas the TTNN errors are generally weaker and more evenly distributed.

The relative $\ell^2$ errors between the learner and $f_4(x)$ in Table \ref{tab:case3.2} agree with these visual observations. With $L=2$, GTNNs already reduce the relative $\ell^2$ error to the $O(10^{-2})$ level under a modest parameter budget, while FNN typically remains larger. When increasing to $L=3$, GTNNs further reach the $O(10^{-3})$ regime, whereas FNN stays around $O(10^{-2})$. Overall, this experiment suggests that GTNNs are more efficient for learning multi-peak Gaussian functions in two dimensions under the same number of parameters as FNN, especially in regions with sharp spatial variation.

\begin{table}[htbp]
\footnotesize
\centering
\caption{Relative $\ell^2$ errors of CTNN, QTNN, TTNN and FNN in Case 3.1.}
\label{tab:case3.1}
\setlength{\tabcolsep}{6pt}
\renewcommand{\arraystretch}{1.15}
\begin{tabular}{c l cccc}
\toprule
$L$ & Parameters & CTNN & QTNN & TTNN & FNN \\
\midrule
\multirow{4}{*}{\textit{$2$}}
 & $\approx100$ & 1.7418e-01 & 2.0220e-01 & 1.7046e-01 & 4.7155e-01 \\ 
 & $\approx300$ & 4.4334e-02 & 4.1000e-02 & 2.7171e-02 & 3.3600e-01 \\ 
 & $\approx400$ & 3.3829e-02 & 4.4792e-02 & 2.7307e-02 & 3.2111e-01 \\ 
 & $\approx500$ & 3.2112e-02 &  \textbf{2.2094e-02} & 2.6662e-02 & 2.7016e-01 \\ 
\midrule
\multirow{4}{*}{\textit{$3$}}
 & $\approx2400$  & 8.3469e-03 & 8.6667e-03 & 8.3353e-03 & 5.8959e-02 \\ 
 & $\approx5400$ & 4.7459e-03 & 4.6993e-03 & 3.9767e-03 & 2.0922e-02 \\ 
 & $\approx9600$ & 4.2536e-03 & 2.6974e-03 & 2.6207e-03 & 1.5424e-02 \\ 
 & $\approx15000$ & 2.4086e-03 & \textbf{2.1760e-03} & 2.2093e-03 & 1.1919e-02 \\ 
\bottomrule
\end{tabular}
\end{table}

\begin{table}[htbp]
\footnotesize
\centering
\caption{Relative $\ell^2$ errors of CTNN, QTNN, TTNN and FNN in Case 3.2.}
\label{tab:case3.2}
\setlength{\tabcolsep}{6pt}
\renewcommand{\arraystretch}{1.15}
\begin{tabular}{c l cccc}
\toprule
$L$ & Parameters & CTNN & QTNN & TTNN & FNN \\
\midrule
\multirow{4}{*}{\textit{$2$}}
 & $\approx100$ & 1.5979e-01 & 1.7382e-01 & 1.6013e-01 & 2.2788e-01 \\ 
 & $\approx400$ & 4.9699e-02 & 6.8795e-02 & 5.9925e-02 & 8.8667e-02 \\ 
 & $\approx600$ & 4.7504e-02 & 5.0366e-02 & 4.9629e-02 & 6.7272e-02 \\ 
 & $\approx700$ & \textbf{3.7997e-02} & 5.0991e-02 & 4.4890e-02 & 5.0608e-02 \\ 
\midrule
\multirow{4}{*}{\textit{$3$}}
 & $\approx4400$  & 1.1563e-02 & 1.3644e-02 & 1.1568e-02 & 2.6661e-02 \\ 
 & $\approx10000$ & 7.5644e-03 & 9.3415e-03 & 8.6212e-03 & 1.7633e-02 \\ 
 & $\approx18000$ & 6.7152e-03 & 6.0824e-03 & 5.5463e-03 & 1.5430e-02 \\ 
 & $\approx28000$ & 5.0771e-03 & 4.9853e-03 & \textbf{4.6452e-03} & 1.3402e-02 \\ 
\bottomrule
\end{tabular}
\end{table}

\subsection{Case 4: Learning the shock waves}
To test the ability of GTNNs in learning functions with sharp transitions, we consider two shock-wave type targets built from a train of smoothed pulses.
We first introduce the sigmoid gate and the $C^{2}$ smoothstep function
\begin{equation}
S_{\varepsilon}(z)=\frac{1}{1+e^{-z/\varepsilon}},\quad 
\eta(u)=6u^{5}-15u^{4}+10u^{3},\quad u\in[0,1].
\end{equation}
Given pulse parameters $(\alpha,a,T_p,T_d)$, define $b=a+T_p$ and $e=b+T_d$, and set
\begin{equation}
u(t)=\min\left\{1,\max\left\{0,\frac{t-b}{T_d}\right\}\right\}.
\end{equation}
We define a 1-D pulse by
\begin{equation}
P(t;\alpha,a,T_p,T_d)
=\alpha\left(S_{\varepsilon_s}(t-a)\left(1-S_{\varepsilon_g}(t-b)\right)
+S_{\varepsilon_g}(t-b)\left(1-S_{\varepsilon_g}(t-e)\right)\left(1-\eta(u(t))\right)\right).
\end{equation}
For the 2-D case, we use the same turn-on/turn-off structure, but replace the smoothstep decay by a steeper sigmoid drop:
\begin{equation}
\widetilde P(t;\alpha,a,T_p,T_d)
=\alpha\left(S_{\varepsilon_s}(t-a)\left(1-S_{\varepsilon_g}(t-b)\right) +S_{\varepsilon_g}(t-b)\left(1-S_{\varepsilon_g}(t-e)\right)\left(1-S_{\varepsilon_f}\left(t-(b+\tau T_d)\right)\right)\right),
\end{equation}
where $\tau\in(0,1)$ controls the drop location inside the decay stage.

The target functions are defined as
\begin{align*}
&\text{Case 4.1: 1-D shock-wave}\quad 
f_5(x)=\sum_{i=1}^{K} P\left(x;\alpha_i,a_i,T_p,T_d\right),\quad x\in[-1,1],\\
&\text{Case 4.2: 2-D shock-wave}\quad 
f_6(x_1,x_2)=A(x_2)\left(\sum_{i=1}^{K} \widetilde P\left(x_1;\alpha_i,a_i,T_p,T_d\right)\right)\,R(x_1),
\quad (x_1,x_2)\in[-1,1]^2,
\end{align*}
where $A(x_2)=\beta_1+\beta_2\exp\left(-(x_2/\sigma_y)^2\right)$ is a smooth envelope in the $x_2$-direction, and $ R(x_1) = 1- S_{\varepsilon_g}\left(x_1-(1-\delta)\right)$ is a gate to ensure $f_6(x_1,x_2)\approx 0$ as $x_1\to 1$. Here $\delta>0$ is a small boundary buffer.

In Case 4, we take $K=10$ pulses, place the starting locations $\{a_i\}$ inside the domain with buffer $\delta$ and a non-overlap constraint, and set the amplitudes to decay geometrically as $\alpha_i=\alpha_0\rho^{\,i-1}$ with $\rho=0.8$. The starting locations $\{a_i\}_{i=1}^K$ are chosen deterministically so that all pulses fit in $[-1,1]$ and do not overlap. For Case 4.1 we use $(T_p,T_d)=(0.08,\,0.045)$ and $(\varepsilon_s,\varepsilon_g)=(5\times10^{-4},\,3\times10^{-3})$. For Case 4.2 we use $(T_p,T_d)=(0.08,\,0.06)$, $(\varepsilon_s,\varepsilon_g,\varepsilon_f)=(10^{-3},\,4\times10^{-3},\,3\times10^{-3})$, and $\tau=0.20$, together with $(\beta_1,\beta_2,\sigma_y)=(0.2,\,0.8,\,0.5)$; we take $\delta=\max\{0.02,\,12\varepsilon_g\}$ so that the right-end gate becomes effective before reaching the boundary. Here $R(x_1)$ is a smooth right-end gate that forces the target to vanish near $x_1=1$. In each case, we select 20000 training feature points $\{x_n\}$ that are evenly distributed within the domain, and fit the dataset $\{x_n,f_k(x_n)\}$ using CTNNs, QTNNs, TTNNs, and FNNs via the least-squares regression~\eqref{14}.

\begin{figure}[htbp]
\centering
\includegraphics[scale=0.6]{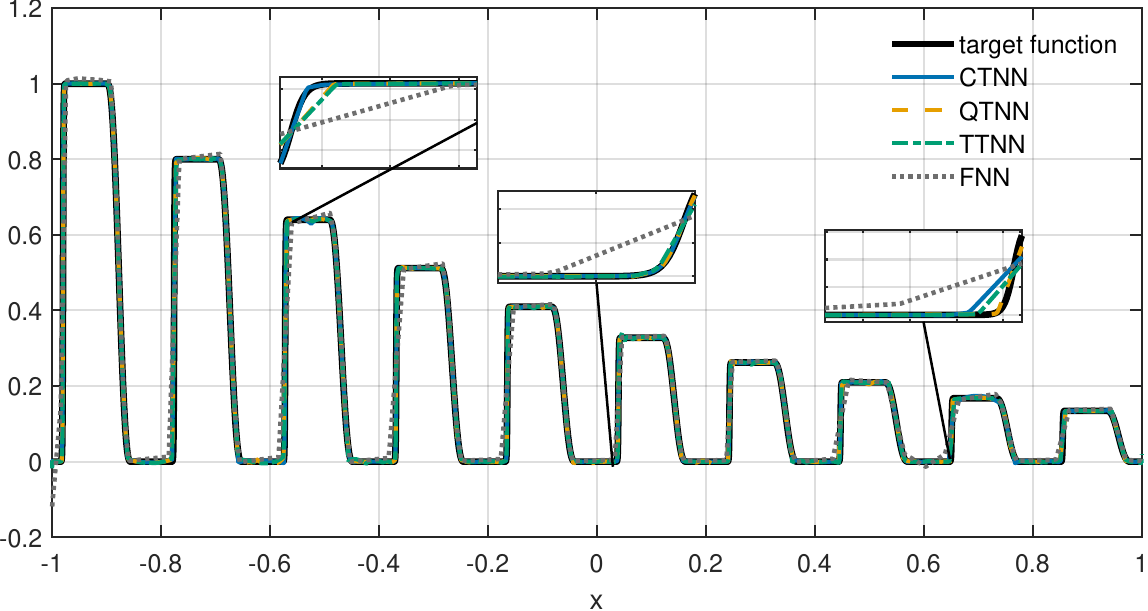}
\caption{\em The target function, the trained FNN, CTNN, QTNN and TTNN in Case 4.1.}
\label{fig:case4.1}
\end{figure}

The 1-D shock-wave $f_5(x)$ is a pulse train with sharp transitions and long flat plateaus, whose amplitudes decay across the interval. As shown in Figure \ref{fig:case4.1}, the GTNN variants track the target waveform closely over the entire domain, and the agreement remains clear in the zoomed-in panels where the jump locations and short transition regions are easiest to inspect. The error trends in Table \ref{tab:case4.1} align with these observations: for $L=2$ all methods stay at the $O(10^{-1})$ level, but GTNNs are consistently more accurate than FNN, and for $L=3$ the GTNN errors improve to $O(10^{-2})$ whereas the FNN remains at $O(10^{-1})$ under the same parameter budget.

\begin{figure}[tbh!]
\centering
\subfloat[target function $f_6(x)$]{
\includegraphics[scale=0.45]{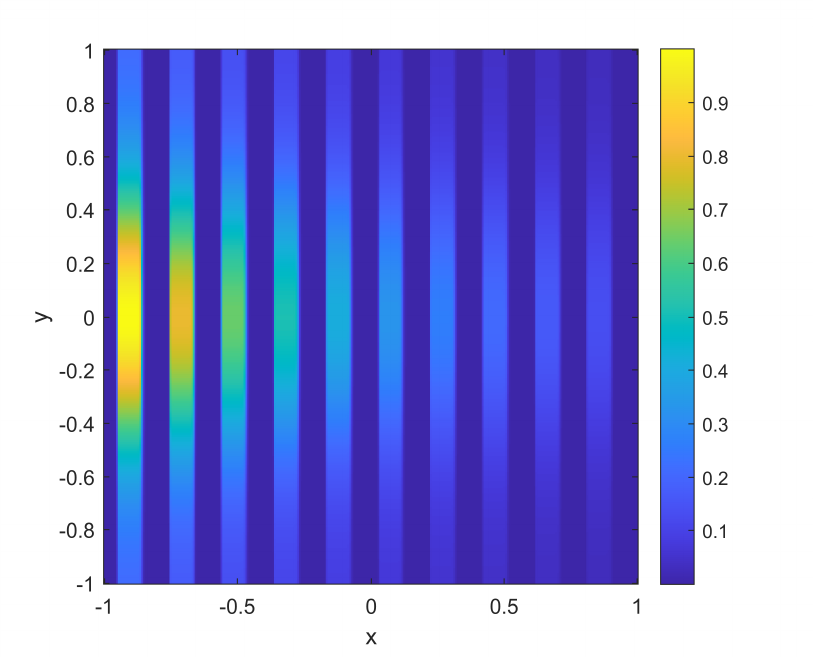}}\hspace{10pt}
\subfloat[$|\phi_{FNN}-f_6(x)|$]{
\includegraphics[scale=0.45]{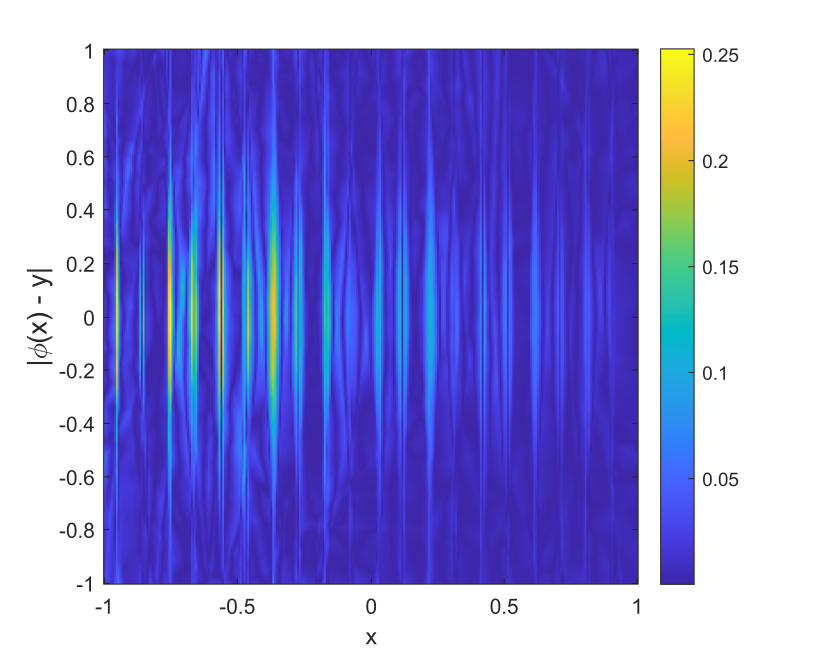}}\hspace{10pt}
\subfloat[$|\phi_{TTNN}-f_6(x)|$]{
\includegraphics[scale=0.45]{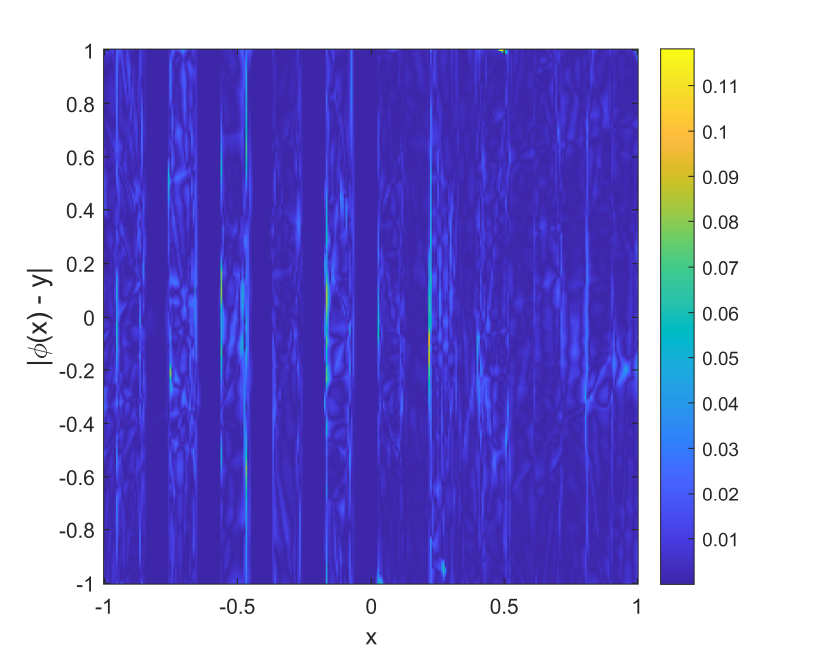}}\hspace{10pt}
\subfloat[$|\phi_{CTNN}-f_6(x)|$]{
\includegraphics[scale=0.45]{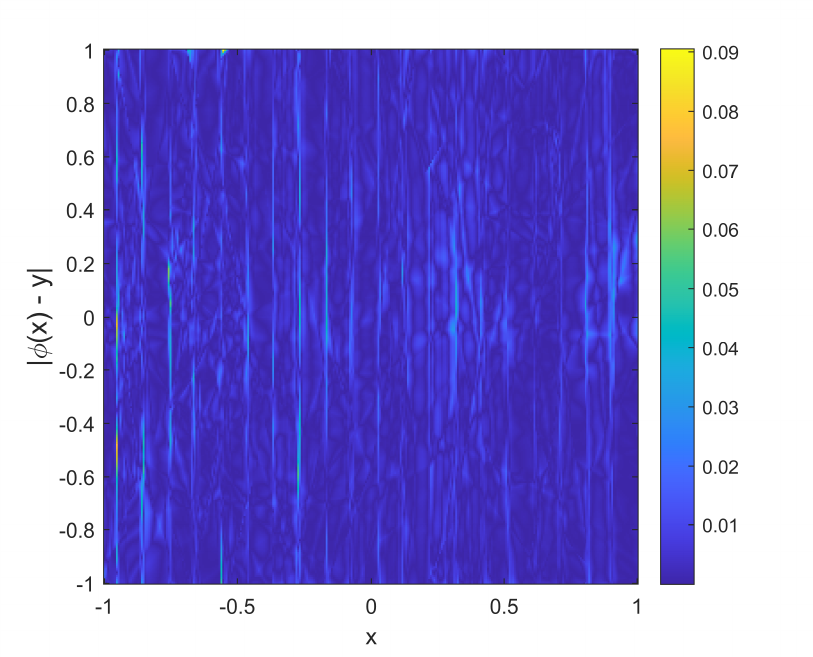}}\hspace{10pt}
\subfloat[$|\phi_{QTNN}-f_6(x)|$]{
\includegraphics[scale=0.45]{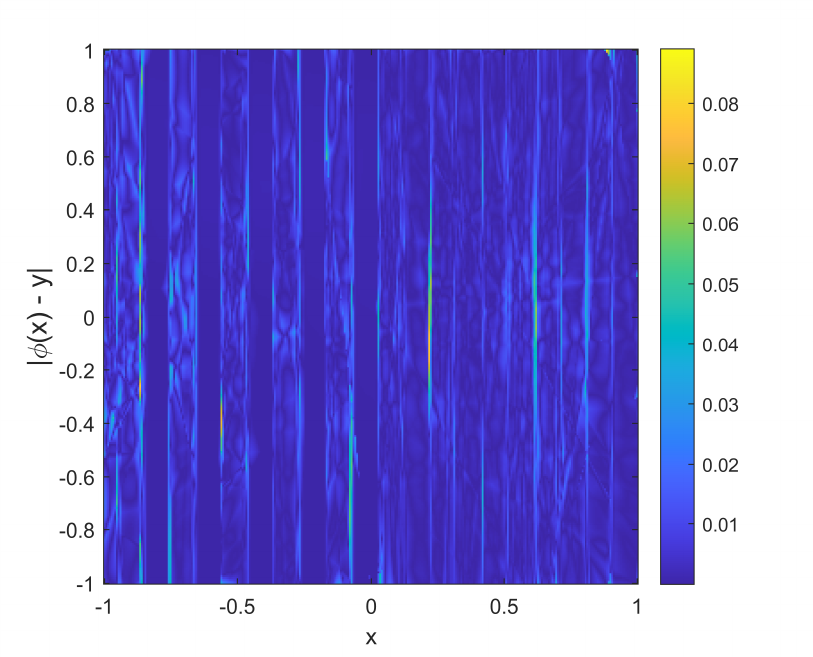}}
\caption{\em 2-D profiles of the target function and absolute errors of the trained FNN, QTNN, CTNN and TTNN in Case 4.2.}
\label{fig:case4.2}
\end{figure}

Case 4.2 extends the same pulse-train mechanism to two dimensions, producing a strongly anisotropic pattern: steep transitions occur primarily along the $x_1$-direction, while the magnitude varies smoothly in $x_2$ through an envelope. This setting tests whether a model can simultaneously capture sharp, aligned interfaces and a smooth modulation without introducing spurious oscillations or excessive diffusion. Figure \ref{fig:case4.2} reports the target profile and pointwise absolute error, where the dominant errors concentrate along the vertical jump lines exactly where the function changes most rapidly. Compared with the FNN, the CTNNs/QTNNs/TTNNs exhibit weaker and more localized error ridges, indicating a better preservation of the sharp fronts and the overall amplitude distribution. As shown in Table \ref{tab:case4.2}, with $L=2$, GTNNs reach the $O(10^{-1})$ regime with smaller constants than FNN, and with $L=3$ the GTNN errors drop to $O(10^{-2})$ while the FNN remains at $O(10^{-1})$. Overall, the 2-D experiment reinforces that, under the same parameter budget, GTNNs achieve higher accuracy on shock-wave targets.

\begin{table}[htbp]
\footnotesize
\centering
\caption{Relative $\ell^2$ errors of CTNN, QTNN, TTNN and FNN in Case 4.1.}
\label{tab:case4.1}
\setlength{\tabcolsep}{6pt}
\renewcommand{\arraystretch}{1.15}
\begin{tabular}{c l cccc}
\toprule
$L$ & Parameters & CTNN & QTNN & TTNN & FNN \\
\midrule
\multirow{4}{*}{\textit{$L=2$}}
 & $\approx100$ & 2.6747e-01 & 3.8263e-01 & 4.3539e-01 & 5.6741e-01 \\ 
 & $\approx300$ & 1.7256e-01 & 2.0752e-01 & 2.1375e-01 & 3.7837e-01 \\ 
 & $\approx400$ & 1.6618e-01 & 1.9415e-01 & 1.6559e-01 & 3.7750e-01 \\ 
 & $\approx500$ & 1.3688e-01 & 1.7997e-01 & \textbf{1.5255e-01} & 3.3728e-01 \\ 
\midrule
\multirow{4}{*}{\textit{$L=3$}}
 & $\approx2000$  & 2.7265e-02 & 6.9668e-02 & 7.6461e-02 & 2.3266e-01 \\ 
 & $\approx18000$ & 3.1518e-02 & 3.4796e-02 & 3.0622e-02 & 1.0045e-01 \\ 
 & $\approx32000$ & 2.3114e-02 &  \textbf{2.0333e-02} & 2.4170e-02 & 8.8366e-02 \\ 
 & $\approx50000$ & 2.2013e-02 & 2.8325e-02 & 3.1589e-02 & 1.0108e-01 \\ 
\bottomrule
\end{tabular}
\end{table}

\begin{table}[htbp]
\footnotesize
\centering
\caption{Relative $\ell^2$ errors of CTNN, QTNN, TTNN and FNN in Case 4.2.}
\label{tab:case4.2}
\setlength{\tabcolsep}{6pt}
\renewcommand{\arraystretch}{1.15}
\begin{tabular}{c l cccc}
\toprule
$L$ & Parameters & CTNN & QTNN & TTNN & FNN \\
\midrule
\multirow{4}{*}{\textit{$L=2$}}
 & $\approx100$ & 3.8428e-01 & 5.1210e-01 & 4.7631e-01 & 6.4765e-01 \\ 
 & $\approx300$ & 2.9716e-01 & 3.2202e-01 & 3.0830e-01 & 5.7002e-01 \\ 
 & $\approx400$ & 2.7876e-01 & 2.9813e-01 & 2.7315e-01 & 5.6313e-01 \\ 
 & $\approx500$ & \textbf{2.5456e-01} & 2.6584e-01 & 2.7676e-01 & 5.3970e-01 \\ 
\midrule
\multirow{4}{*}{\textit{$L=3$}}
 & $\approx2000$  & 5.9547e-02 & 1.4730e-01 & 7.6317e-02 & 2.8788e-01 \\ 
 & $\approx18000$ & 4.8498e-02 & 7.2987e-02 & 8.2592e-02 & 2.0016e-01 \\ 
 & $\approx32000$ & 4.1670e-02 & 4.5155e-02 & 6.7626e-02 & 1.7071e-01 \\ 
 & $\approx50000$ & \textbf{2.9686e-02} & 3.4528e-02 & 3.6917e-02 & 1.5966e-01 \\ 
\bottomrule
\end{tabular}
\end{table}

\section{Conclusion}
We develop a novel class of DNN-like parametrized functions, i.e., general transformation neural networks (GTNNs), for high-dimensional approximation. The main idea is to generalize the affine transformation $x\mapsto Wx+b$ of conventional abstract neurons to more general transformations $x\mapsto (Wx+b)\psi(x;\htheta)$, where $\psi$ is a specific parametrized function. For ReLU activation, it is equivalent to generalizing the linear shape functions of the conventional DNNs to more complex shape functions. We take $\psi(x;\htheta)$ to be polynomials of degree 2 and 3, and the finite Fourier series for a detailed study, establishing universal approximation theorems and error bound estimation. In numerical experiments, we consider least squares regression problems, where GTNNs show significant advantages over fully connected neural networks; higher accuracy is obtained in learning oscillatory functions, the flow of the shock-turbulence problem, multi-peak Gaussian functions and shock waves. 

Note that GTNNs are proposed based on the idea of the mesh-based approximations, such as finite elements and spectral elements. To implement the meshing, we use the ReLU activation, which vanishes in one part of the domain and is the identity in the other. More generally, one can adopt other $L$-shaped activation functions such as the softplus function and the exponential linear unit. Whereas it is questionable whether sigmoidal functions are as effective as $L$-shaped activation functions, and this could be investigated in future work.

Furthermore, in this work, we consider only three simple types of GTNNs: QTNNs, CTNNs and TTNNs, which use quadratic polynomials, cubic polynomials and finite Fourier series as the activated transformation, respectively. We think the applicable GTNN types extend well beyond these three. It might be promising to investigate other types of general transformations, such as wavelet functions or Gaussian functions, for specific problems.

\vspace{10pt}
\noindent{\fontsize{9.0pt}{\baselineskip}\selectfont \textbf{Funding} This work is supported by National Natural Science Foundation of China Major Research Plan (92370101).}

\vspace{10pt}
\noindent{\fontsize{9.0pt}{\baselineskip}\selectfont\textbf{Data availability} Data will be made available on request.}

\vspace{10pt}
\noindent{\fontsize{9.0pt}{\baselineskip}\selectfont\textbf{Declaration of competing interest} The authors declare that they have no known competing financial interests or personal relationships that could have appeared to influence the work reported in this paper.}

\appendix 

\section{Proofs for Lemmas and Theorems}\label{proof}
\subsection{Proof of Lemma \ref{lem1}}
\label{prflem1}
We consider $\hat{f}$ to be the Fourier transform of $f$. By the definition \eqref{18}, we have
\begin{equation}\label{19}
\Gamma(f)= \int_{\mbR^d} \|\omega\|_2^2|\hat{f}(\omega)|\td \omega.
\end{equation}
By using the truth $\|\omega\|_2^2|\hat{f}(\omega)| = \frac{1}{\sqrt{1+\|\omega\|_2^{2s-4}}}\sqrt{1+\|\omega\|_2^{2s-4}}\|\omega\|_2^2|\hat{f}(\omega)|$ and Cauchy–Schwarz inequality, we have
\begin{equation}\label{20}
\int_{\mbR^d} \|\omega\|_2^2|\hat{f}(\omega)|\td \omega \leq 
\left( \int_{\mbR^d}\frac{\td \omega}{1+\|\omega\|_2^{2s-4}}\right)^{\frac{1}{2}}\left(\int_{\mbR^d}(1+\|\omega\|_2^{2s-4}) \|\omega\|_2^{4} |\hat{f}(\omega)|^2 \td \omega\right)^{\frac{1}{2}}.
\end{equation}

Because $2s-4>d$, we have $ \int_{\mbR^d}\frac{\td \omega}{1+\|\omega\|_2^{2s-4}}< \infty$. We denote
\begin{equation}
I_1 = \int_{\mbR^d}\|\omega\|_2^{4}|\hat{f}(\omega)|^2 \td \omega, ~ I_2 = \int_{\mbR^d}\|\omega\|_2^{2s}|\hat{f}(\omega)|^2 \td \omega.
\end{equation} 
By the derivative properties of the Fourier transform and the Plancherel theorem, we have
\begin{multline}\label{30}
I_1 = \int_{\mbR^d}\|\omega\|_2^{4}|\hat{f}(\omega)|^2 \td \omega = \int_{\mbR^d}(\omega_1^2 + \omega_2^2 + \cdots + \omega_d^2)^2|\hat{f}(\omega)|^2 \td \omega =  \\\int_{\mbR^d} \left(\sum_{i,j=1}^d \omega_i^2\omega_j^2\right)|\hat{f}(\omega)|^2 \td \omega 
= \int_{\mbR^d} \sum_{i,j=1}^d|(i\omega_i)(i\omega_j)\hat{f}(\omega)|^2 \td \omega = \\\int_{\mbR^d} \sum_{i,j=1}^d|\widehat{\partial_{x_i,x_j}f(\omega)}|^2\td \omega = (2\pi)^d\sum\limits_{\|\alpha\|_1=2}\|D^\alpha f\|_{L^2}^2.
\end{multline}

For $I_2$, by a similar argument, we have
\begin{equation}\label{31}
I_2 = (2\pi)^d\sum\limits_{\|\beta\|_1=s} \binom{s}{\beta}\|D^\beta f\|_{L^2}^2 \leq \max\limits_{\|\beta\|_1=s}\binom{s}{\beta} (2\pi)^d\sum\limits_{\|\beta\|_1=s}\|D^\beta f\|_{L^2}^2
\end{equation}
where $\binom{s}{\beta} = \frac{s!}{\beta_1!\dots\beta_d!}$ denotes the multinomial coefficient corresponding to the multi–index $\beta$. Combining \eqref{30} and \eqref{31}, we have
\begin{equation}
\Gamma(f) \leq C(d,s)\left(\sum\limits_{\|\alpha\|_1=2}\|D^\alpha f\|_{L^2}^2 + \sum\limits_{\|\beta\|_1=s}\|D^\beta f\|_{L^2}^2\right)^\frac{1}{2},
\end{equation}
where $C(d,s)= (2\pi)^{d/2} \left( \int_{\mbR^d}\frac{\td \omega}{1+\|\omega\|_2^{2s-4}}\right)^{\frac{1}{2}} \max\limits_{\|\beta\|_1=s}\binom{s}{\beta}$.

\subsection{Proof of Lemma \ref{lem2}}\label{prflem2}
For any $x\in\Omega$ and $f\in H^s(\mbR^d)$, we have the following fact
\begin{equation}\label{36}
    f(x) - f(0) - \nabla f^\top(0) x = \operatorname{Re} \int_{\mbR^d} (e^{i\omega \cdot x} - i\omega \cdot x -1)\hat{f}(\omega)\td \omega.
\end{equation}   
For any $c\in\mbR^+$, if $|z|\leq c$, we note the identity
\begin{equation}\label{37}
    -\int_{0}^{c}\sigma(z-u) e^{iu} + \sigma(-z-u) e^{-iu} \td u = e^{iz} - iz - 1.
\end{equation}
For all $\omega\in \mbR^d \backslash \{0\}$, $x\in \mbR^d$, we let $z = \omega^\top x$, $a = \omega/\|\omega\|_2$ and $u = \|\omega\|_2t$, $t\in[0,1]$ and substitute them into \eqref{36}. Then we have
\begin{equation}\label{38}
    -\|\omega\|_{2}^{2} \int_{0}^{1} \sigma(a \cdot x - t) e^{i\|\omega\|_{2} t} + \sigma(-a \cdot x - t) e^{-i\|\omega\|_{2} t}  \td t = e^{i \omega \cdot x} - i \omega \cdot x - 1.
\end{equation}
We substitute \eqref{38} into \eqref{36} using Fubini’s theorem and have
\begin{multline}\label{39}
    f(x) - f(0) - \nabla f^\top(0) x = \\
    -\operatorname{Re} \int_{\mbR^d} \int_{0}^{1} \left( \sigma(a \cdot x - t) e^{i\|\omega\|_{2} t} + \sigma(-a \cdot x - t) e^{-i\|\omega\|_{2} t} \right) \|\omega\|_2^2 \hat{f}(\omega) \td t\td \omega. 
\end{multline} 

By Lemma \ref{lem1}, for any $f\in H^s(\mbR^d)$ with $s = \lfloor d/2\rfloor + 3$, we have $\Gamma(f)<\infty$. Therefore, the above integral \eqref{39} is well-defined. Let $\hat{f}(\omega) = e^{i\gamma(\omega)}|\hat{f}(\omega)|$, where $\gamma(w)$ means the argument of $\hat{f}(\omega)$. Taking the real part, \eqref{39} becomes
\begin{equation}\label{40}
   -\sum_{z=\pm 1}\int_{\mbR^d} \int_{0}^{1}  \sigma(za \cdot x - t)\cos(z\|\omega\|_2 t + \gamma(\omega)) \|\omega\|_2^2 \hat{f}(\omega) \td t\td \omega. 
\end{equation} 

Moreover, we can consider a density function on $\{-1,1\} \times [0,1] \times \mbR^d \backslash \{0\}$ defined by
\begin{equation}\label{41}
   p(z,t,\omega) = \frac{\left| \cos\left( z \|\omega\|_{2} t + \gamma(\omega) \right) \right| \|\omega\|_{2}^{2} |\hat{f}(\omega)|}{v},
\end{equation} 
where 
\begin{multline}
v = \int_{\mathbb{R}^{d}} \int_{0}^{1} \left( \left| \cos \left( \|\omega\|_{2} t + \gamma(\omega) \right) \right| + \left| \cos \left( \|\omega\|_{2} t - \gamma(\omega) \right) \right| \right) \|\omega\|_{2}^{2} |\hat{f}(\omega)|  \td t  \td\omega\\
\leq 2\Gamma(f)<\infty.
\end{multline} 
The corresponding probability measure is denoted as $P$. We define a sign function $\eta(z,t,\omega) = -\operatorname{sgn}(\cos\left( z \|\omega\|_{2} t + \gamma(\omega) \right))$ and combining \eqref{40} and the definition of $p(z,t,\omega)$, \eqref{39} can be written as
\begin{equation}\label{42}
   f(x) - f(0) - \nabla f^\top(0) x = v\mathbb{E}_{(z,t,\omega)\sim P}[\eta(z,t,\omega)\sigma(za(\omega) \cdot x - t)],
\end{equation} 

Because for $\htheta\in \Theta$, $\psi(\cdot;\htheta) \geq0$ and $\sigma(x)=\mathrm{ReLU}(x)$, $\psi$ satisfies $\sigma(\psi t) = \psi \sigma(t)$. For \eqref{42}, by Fubini’s theorem and $\int_{\Theta} \psi(x;\htheta) \nu\td (\htheta) =1$, we have
\begin{multline}
f(x) - f(0) - \nabla f^\top(0) x = v\mathbb{E}_{(z,t,\omega)\sim P}\left[\int_{\Theta} \psi(x;\htheta) \nu\td (\htheta)~\eta(z,t,\omega)\sigma(za(\omega) \cdot x - t)\right]\\
=v\mathbb{E}_{(z,t,\omega,\htheta)\sim (P \times \nu)}\left[\eta(z,t,\omega)\sigma\left((za(\omega) \cdot x - t)\psi(x;\htheta)\right)\right].
\end{multline} 
By the Radon-Nikodym theorem, we can define the following signed measure on $\mathcal{D}_1 = \{-1,1\} \times [0,1] \times \mbR^d \times \Theta$
\begin{equation}\label{43}
   \mu = v\eta\cdot(P \times \nu). 
\end{equation} 
On the other hand, the mapping
\begin{equation}\label{44}
    \Phi: \mathcal{D}_1\to\mathcal{D}_0,\quad (z,t,\omega,\htheta)\to (w,b,\htheta) = (za(\omega),-t,\htheta)
\end{equation}
is measurable. Then, we can take $\rho^* :=\mu\circ\Phi^{-1}$ as the pushforward measure of $\mu$, i.e., $\forall A\subset\mathcal{D}_0$ is measurable then $\rho^*(A) := \mu\left( \Phi^{-1}(A) \right)$. Moreover, we have
\begin{multline}\label{45}
\int_{\mathcal{D}_0} \sigma\left((w^\top x+b) \psi(x;\htheta)\right) \rho^*(\td w,\td b, \td \htheta) = \\  
\int_{\mathcal{D}_1} \sigma\left((za(\omega) \cdot x - t) \psi(x;\htheta)\right) \mu(\td z, \td \omega,\td t, \td \htheta) = f(x) - f(0) - \nabla f^\top(0)x. 
\end{multline}

For any finite signed measure $\zeta$ on $\mathcal{D}$, we define the total variation norm  $\|\zeta\|_{\operatorname{TV}}=|\zeta|(\mathcal{D})$, where $|\zeta|$ is total variation measure of $\zeta$. Because $|\eta|=1$ and $P \times \nu$ is joint probability measure, we have $\|\mu\|_{\operatorname{TV}}=v\leq 2\Gamma(f)$. Furthermore, we have $\| \rho^* \|_{\operatorname{TV}} \leq \|\mu\|_{\operatorname{TV}}\leq 2\Gamma(f)$. We consider the linear term $-\nabla f^\top(0)x$ in \eqref{45}. For $j=1,\dots,d$, we let $e_j=(0,\dots,1,\dots,0)$ be the one-hot vector whose $j$-th component is one and the others are zero. By using the identity for $\sigma=\mathrm{ReLU}$, $u = \sigma(u)-\sigma(-u)$, we have
\begin{equation}\label{46}
     -\nabla f^\top(0) x = -\sum_{j=1}^d \partial_jf(0)(\sigma(e_j^\top x) -\sigma(-e_j^\top x)).
\end{equation}
If we let $\delta_{(w,b)}$ be the Dirac measure whose support point is $(w,b)$ on $\mbR^d \times \mbR$. Then, we let the joint probability measure 
\begin{equation}\label{47}
     \rho_1 = (-\sum_{j=1}^d\partial_jf(0)(\delta_{(e_j,0)} -\delta_{(-e_j,0)}))\times \nu.
\end{equation}

For any finite signed measure $\zeta_1$ and $\zeta_2$, we have the fact
\begin{equation}
    |\zeta_1\times \zeta_2|(E\times F) = |\zeta_1|(E) \cdot|\zeta_2|(F),
\end{equation}
where $E$ and $F$ are measurable sets in the domains of $\zeta_1$ and $\zeta_2$, respectively. Using the fact, we have
\begin{multline}
    \|\rho_1\|_{\operatorname{TV}}=|\rho_1|(\mbR^d \times \mbR \times \Theta) = \left| (-\sum_{j=1}^d\partial_jf(0)(\delta_{(e_j,0)} -\delta_{(-e_j,0)}))\times \nu\right|(\mbR^d \times \mbR \times \Theta)\\
    =\left|\sum_{j=1}^d\partial_jf(0)(\delta_{(e_j,0)} -\delta_{(-e_j,0)})\right|(\mbR^d \times \mbR) \cdot |\nu|(\Theta)\\
    =\int_{\mbR^d \times \mbR}\left|\sum_{j=1}^d\partial_jf(0)(\delta(x-(e_j,0)^\top) -\delta(x+(e_j,0)^\top))\right| \td x \cdot1=2\sum_{j=1}^d|\partial_j f(0)| =  2\|\nabla f(0)\|_1
\end{multline}
Now we can represent the linear term in the following integral form
\begin{equation}\label{48}
     -\nabla f(0)^\top x =\int_{\mathcal{D}_0} \sigma\left((w^\top x+b) \psi(x;\htheta)\right) \rho_1(\td w,\td b, \td \htheta). 
\end{equation}

For the constant term $-f(0)$ in \eqref{45}, we can similarly define the probability measure $\rho_2 = -f(0) \delta_{(0,1)} \times \nu$ such that 
\begin{equation}\label{49}
     -f(0) = \int_{\mathcal{D}_0} \sigma\left((w^\top x+b) \psi(x;\htheta)\right) \rho_2(\td w,\td b, \td \htheta)
\end{equation}
with $\|\rho_2\|_{\operatorname{TV}} = |f(0)|$. We let $\rho = \rho^* + \rho_2 + \rho_1$ and combining \eqref{45}, \eqref{48} and \eqref{49} we obtain 
\begin{equation}\label{50}
    f(x) = \int_{\mathcal{D}_0} \sigma((w^\top x+b)\psi(x;\htheta)) \rho(\td w, \td b, \td \htheta),\quad \forall x\in \Omega,
\end{equation}
where $\rho$ satisfies
\begin{equation}
    \|\rho\|_{\operatorname{TV}} \leq \|\rho^*\|_{\operatorname{TV}} + \|\rho_2\|_{\operatorname{TV}} + \|\rho_1\|_{\operatorname{TV}}\leq 2\Gamma(f) + 2\|\nabla f(0)\|_1 + |f(0)|.
\end{equation}
\subsection{Proof of Theorem \ref{thm3}}\label{prfthm3}
By Lemma \ref{lem1} we know that for $s = \lfloor d/2\rfloor + 3$ there exists a constant $C_1=C_1(d,s)$ such that
\begin{equation}\label{55}
  \Gamma(f) \le C_1\,\|f\|_{H^s(\mbR^d)}.
\end{equation}
Moreover, since $s>d/2+1$, Sobolev embedding yields $H^s(\mbR^d)\hookrightarrow C_{\mathrm{b}}^1(\mbR^d)$, where $C_{\mathrm{b}}^1(\mbR^d)$ denotes the space of continuously differentiable functions $f$: $\mbR^d \rightarrow \mbR$ whose values and all first–order partial derivatives are bounded. Then there exists $C_2=C_2(d,s)$ such that
\begin{equation}\label{56}
  |f(0)| + \|\nabla f(0)\|_1 \le C_2\,\|f\|_{H^s(\mbR^d)}. 
\end{equation}
Combining \eqref{55} and \eqref{56} we obtain
\begin{equation}\label{63}
  2\Gamma(f) + 2\|\nabla f(0)\|_1 + |f(0)| \le C\,\|f\|_{H^s(\mbR^d)}
\end{equation}
for some constant $C=C(d,s)$.

By Lemma \ref{lem2} there exists a finite signed measure $\rho$ on $\mathcal{D}_0 = \mbR^d\times\mbR\times\Theta$ such that
\begin{equation}\label{eq:f-rho}
  f(x) = \int_{\mathcal{D}_0} \sigma\big((w^\top x + b)\psi(x;\htheta)\big)\, \rho(\td w,\td b,\td \htheta),~\forall x\in\Omega,
\end{equation}
and $\|\rho\|_{\mathrm{TV}} \le 2\Gamma(f) + 2\|\nabla f(0)\|_1 + |f(0)|$. Moreover, by the explicit construction in the proof of Lemma \ref{lem2}, the support of $\rho$ is contained in the set $\{(w,b,\htheta)\in \mathcal{D}_0: |w|\le 1,|b|\le 1\}$.

By the Radon-Nikodym theorem, there exists a measurable function $p$ such that 
\begin{equation}
    \td \rho = p\td|\rho|.
\end{equation}
where $p:\mathcal{D}_0\to\{-1,1\}$ is a measurable sign function given by the polar decomposition of $\rho$ with respect to $|\rho|$. Define the probability measure 
\begin{equation}
    \mu_0 := \frac{1}{\|\rho\|_{\mathrm{TV}}}|\rho|
\end{equation}
on $\mathcal{D}_0$, so that $\mu_0(\mathcal{D}_0)=1$. On the enlarged parameter space $\mathcal{D}= \mbR\times\mbR^d\times\mbR\times\Theta$ we define a measurable map
\begin{equation}
  \Psi:\mathcal{D}_0\to \mathcal{D},\quad \Psi(w,b,\htheta) = (a,w,b,\htheta),
\end{equation}
where $a=\|\rho\|_{\mathrm{TV}}\,p(w,b,\htheta)$. Let $\pi$ be the pushforward measure, i.e., $\pi := \mu_0 \circ\Psi^{-1}$. Then $\pi$ is a probability measure on $\mathcal{D}$ and for any bounded measurable $G:D\to\mbR$ we have
\begin{equation}\label{57}
\int_\mathcal{D} G(a,w,b,\htheta)\pi(\td a,\td w,\td b,\td \htheta) = \int_{\mathcal{D}_0} G\big(\|\rho\|_{\mathrm{TV}}p(w,b,\htheta),w,b,\htheta\big)\mu_0(\td w,\td b,\td \htheta).
\end{equation}

In particular, taking 
\begin{equation}
 G(a,w,b,\htheta) = a\sigma\big((w^\top x + b)\psi(x;\htheta)\big),
\end{equation}
and using $d\rho = pd|\rho| = \|\rho\|_{\mathrm{TV}}pd\mu_0$, we obtain
\begin{multline}
\int_\mathcal{D} a\sigma\big((w^\top x + b)\psi(x;\htheta)\big)\pi(\td a,\td w,\td b,\td \htheta) \\
= \int_{\mathcal{D}_0} \|\rho\|_{\mathrm{TV}}p(w,b,\hat\theta)\sigma\big((w^\top x + b)\psi(x;\hat\theta)\big)\mu_0(\td w,\td b,\td \htheta)\\
= \int_{\mathcal{D}_0} \sigma\big((w^\top x + b)\psi(x;\htheta)\big)\rho(\td w,\td b,\td \htheta) = f(x).
\end{multline}
Thus $\pi\in\Pi_f$ in the sense of \eqref{B_norm_GTNN}. Since $(w,b,\htheta)$ are unchanged under $\Psi$, the support of $\pi$ also satisfies $|w|\le 1$ and $|b|\le 1$.

Using \eqref{57} with
\begin{equation}
G(a,w,b,\htheta)= a^2 M_\psi(\htheta)^2(|w|+|b|)^2,
\end{equation}
and recalling that $a = \|\rho\|_{\mathrm{TV}}p(w,b,\hat\theta)$ and $|p|=1$, we obtain
\begin{multline}
I(f):= \int_\mathcal{D} a^2 M_\psi(\htheta)^2(|w|+|b|)^2\pi(\td a,\td w,\td b,\td \htheta)\\
 = \int_{\mathcal{D}_0} \|\rho\|_{\mathrm{TV}}^2 M_\psi(\htheta)^2(|w|+|b|)^2\mu_0(\td w,\td b,\td \htheta)\\
=\|\rho\|_{\mathrm{TV}}\int_{\mathcal{D}_0} M_\psi(\htheta)^2(|w|+|b|)^2|\rho|(\td w,\td b,\td \htheta).
\end{multline}
By the property $|w|\le 1$, $|b|\le 1$ on $\mathrm{supp}(|\rho|)$ we have
\begin{equation}
  (|w|+|b|)^2 \le (1+1)^2 = 4.
\end{equation}
Therefore,
\begin{equation}\label{58}
  I(f) \le 4\|\rho\|_{\mathrm{TV}}\int_{\mathcal{D}_0} M_\psi(\hat\theta)^2|\rho|(\td w,\td b,\td \htheta).
\end{equation}

Next, recall from the proof of Lemma \ref{lem2} that $\rho$ is the sum of three parts,
\begin{equation}
  \rho = \rho^* + \rho_1 + \rho_2,
\end{equation}
corresponding to the Fourier part, the linear term $\nabla f(0)^\top x$ and the constant term $f(0)$, respectively. Consequently,
\begin{equation}
|\rho| \le |\rho^*| + |\rho_1| + |\rho_2|, 
\end{equation}
and hence
\begin{equation}
  \int_{\mathcal{D}_0} M_\psi(\htheta)^2|\rho|(\td w,\td b,\td \htheta) \le I_1 + I_2 + I_3,
\end{equation}
where
\begin{align}
I_1 = \int_{\mathcal{D}_0} M_\psi(\htheta)^2|\rho^*|(\td w,\td b,\td \htheta),\\ 
I_2 = \int_{\mathcal{D}_0} M_\psi(\htheta)^2\,|\rho_1|(\td w,\td b,\td \htheta),\\
I_3 = \int_{D_0} M_\psi(\htheta)^2\,|\rho_2|(\td w,\td b,\td \htheta).
\end{align}

By construction of the proof of Lemma \ref{lem2}, we have
\begin{equation}
\rho^* = v\,\eta\,(P\times\nu)\circ\Phi^{-1},
\end{equation}
where $P$ is a probability measure on the auxiliary variables $(z,t,\omega)$, $\eta\in\{-1,1\}$ and $v\le 2\Gamma(f)$. 
In particular, the $\htheta$–marginal of $|\rho^*|$ is exactly $v\nu$, and we obtain
\begin{equation}\label{59}
  I_1 = \int_\Theta M_\psi(\hat\theta)^2\,v\nu(\td \htheta)= vJ_{\psi,\nu} \le 2\Gamma(f)J_{\psi,\nu}.
\end{equation}

By the definition of $\rho_1$ in the proof of Lemma \ref{lem2}, we have
\begin{equation}
    |\rho_1|= \left(\sum_{j=1}^d |\partial_j f(0)|(\delta_{(e_j,0)}+\delta_{(-e_j,0)})\right)\times\nu,
\end{equation}
and therefore
\begin{equation}\label{60}
  I_2= 2\sum_{j=1}^d |\partial_j f(0)|\int_\Theta M_\psi(\htheta)^2\nu(\td \htheta)= 2\|\nabla f(0)\|_1J_{\psi,\nu}.
\end{equation}
For the last term, we have
\begin{equation}
\rho_2 = f(0)\delta_{(0,1)}\times\nu,\quad |\rho_2| = |f(0)|\delta_{(0,1)}\times\nu,
\end{equation}
and thus
\begin{equation}\label{61}
  I_3 = |f(0)|\int_\Theta M_\psi(\htheta)^2\,\nu(d\htheta) = |f(0)|J_{\psi,\nu}.
\end{equation}

Combining \eqref{59},\eqref{60} and \eqref{61} we deduce
\begin{equation}
\int_{\mathcal{D}_0} M_\psi(\htheta)^2|\rho|(\td w,\td b,\td \htheta) \le J_{\psi,\nu}(2\Gamma(f) + 2\|\nabla f(0)\|_1 + |f(0)|).
\end{equation}
Substituting this into \eqref{58}, and using $\|\rho\|_{\mathrm{TV}}\le 2\Gamma(f) + 2\|\nabla f(0)\|_1 + |f(0)|$, we obtain
\begin{equation}\label{62}
  I(f) \le 4J_{\psi,\nu}(2\Gamma(f) + 2\|\nabla f(0)\|_1 + |f(0)|)^2.
\end{equation}
By definition \eqref{B_norm_GTNN} and the fact that $\pi\in\Pi_f$ we have
\begin{equation}
\mathcal{C}_\psi^2(f)= \inf_{\tilde\pi\in\Pi_f}\int_\mathcal{D} a^2 M_\psi(\htheta)^2(|w|+|b|)^2\,\tilde\pi(\td a,\td w,\td b,\td \htheta)\le I(f).
\end{equation}
Combining this with \eqref{62} and \eqref{63}, we obtain
\begin{equation}
\mathcal{C}_\psi^2(f)\le 4J_{\psi,\nu}(2\Gamma(f) + 2\|\nabla f(0)\|_1 + |f(0)|)^2\le 4J_{\psi,\nu}C^2\ \|f\|_{H^s(\mathbb{R}^d)}^2.
\end{equation}
Taking square roots gives
\begin{equation}
\mathcal{C}_\psi(f)\le 2C(d,s)\sqrt{J_{\psi,\nu}}\,\|f\|_{H^s(\mathbb{R}^d)} < \infty, 
\end{equation}
Therefore, $\mathcal{C}_\psi(f)$ is finite and hence $f\in\mathcal{B}$. Since the above bound holds for every
$f\in H^s(\mathbb{R}^d)$, we conclude that
\begin{equation*}
H^s(\mathbb{R}^d)\subset \mathcal{B},
\end{equation*}
where $s=\lfloor d/2\rfloor+3$.

\bibliographystyle{plainnat}
\bibliography{references}

@article{Pinkus1999,
  title={Approximation theory of the {MLP} model in neural networks},
  author={A. Pinkus },
  journal={Acta Numer.},
  volume={8},
  pages={143--195},
  year={1999}
}

@article{Fan2020,
    title = {Universal approximation with quadratic deep networks},
    author = {F. L. Fan and J. Xiong and G. Wang},
    journal = {Neural Netw.},
    volume = {124},
    pages = {383--392},
    year = {2020}
}

@article{Fan2018,
    title = {A new type of neurons for machine learning},
    author = {F. L. Fan and W. Cong and G. Wang},
    journal = {Int. J. Numer. Methods Biomed. Eng.},
    volume = {34},
    pages = {e2920},
    year = {2018}
}

@article{Fan2017,
    title = {Generalized backpropagation algorithm for training second-order neural networks},
    author = {F. L. Fan and W. Cong and G. Wang},
    journal = {Int. J. Numer. Methods Biomed. Eng.},
    volume = {34},
    pages = {e2956},
    year = {2017}
}

@article{Fan2025,
  title={One neuron saved is one neuron earned: on parametric efficiency of quadratic networks},   
  author={F. L. Fan and H. C. Dong and Z. Wu and  L. Ruan and T. Zeng and Y. Cui and J. X. Liao},
  journal={IEEE Trans. Pattern Anal. Mach. Intell.}, 
  volume={47},
  pages={9702--9717},
  year={2025}
}

@article{Fan2023,
  title={On expressivity and trainability of quadratic networks}, 
  author={F. L. Fan and M. Li and F. Wang and R. Lai and G. Wang},
  journal={IEEE Trans. Neural Netw. Learn. Syst.}, 
  volume={36},
  pages={1228--1242},
  year={2025}
}

@inproceedings{He2015,
  title={Delving deep into rectifiers: {S}urpassing human-level performance on imagenet classification},
  author={K. He and X. Zhang and S. Ren and J. Sun},
  booktitle={Proceedings of the IEEE international conference on computer vision},
  year={2015}
}

@inproceedings{Trottier2017,
  title={Parametric exponential linear unit for deep convolutional neural networks},
  author={L. Trottier and P. Giguere and B. Chaib-draa},
  booktitle={IEEE International Conference on Machine Learning and Applications},
  year={2017}
}

@inproceedings{Jin2016,
  title   = {Deep learning with {S}-shaped rectified linear activation units},
  author  = {X. Jin and C. Xu and J. Feng and Y. Wei and  J. Xiong and S. Yan},
  booktitle = {Proceedings of the AAAI Conference on Artificial Intelligence},
  year    = {2016}
}

@article{Lu2021,
  title = {Deep network approximation for smooth functions},
  author = {J. Lu and Z. Shen and H. Yang and S. Zhang},
  journal = {SIAM J. Math. Anal.},
  volume = {53},
  pages = {5465--5506},
  year = {2021}
}

@article{Barron1993,
  title  = {Universal approximation bounds for superpositions of a sigmoidal function},
  author = {A. R. Barron},
  journal = {IEEE Trans. Inform. Theory},
  volume = {39},
  pages = {930--945},
  year  = {1993}
}

@article{Titarev2004,
  title={Finite-volume {WENO} schemes for three-dimensional conservation laws},
  author={V. A. Titarev and E. F. Toro},
  journal={J. Comput. Phys.},
  volume={201},
  pages={238--260},
  year={2004}
}

@article{Hornik1989,
    title = {Multilayer feedforward networks are universal approximators},
    author = {K. Hornik and M. Stinchcombe and H. White},
    journal = {Neural Netw.},
    volume = {2},
    pages = {359--366},
    year = {1989}
}

@article{Cybenko1989,
    title = {Approximation by superpositions of a sigmoidal function},
    author = {G. Cybenko},
    journal = {Math. Control Signals Systems},
    volume = {2},
    pages = {303--314},
    year = {1989}
}

@article{E2019,
  title = {A priori estimates of the population risk for two-layer neural networks},
  author = {W. E and C. Ma and L. Wu},
  journal = {Commun. Math. Sci.},
  volume = {17},
  pages = {1407--1425},
  year = {2019}
}

@article{E2022,
  title={The {B}arron space and the flow-induced function spaces for neural network models},
  author={W. E and C. Ma and L. Wu},
  journal={Constr. Approx.},
  volume={55},
  pages={369--406},
  year={2022}
}

@article{E2018,
    title  = {Exponential convergence of the deep neural network approximation for analytic functions},
    author = {W. E and Q. Wang},
    journal = {Sci. China Math.},
    volume  = {61},
    pages = {1733--1740},
    year = {2018}
}

@article{E2020,
   title = {Towards a mathematical understanding of neural network-based machine learning: what we know and what we don't}, 
   author = {W. E and C. Ma and L. Wu and S. Wojtowytsch},
   journal = {CSIAM Trans. Appl. Math.},
   volume = {1},
   pages = {561--614},
   year = {2020}
}

@inproceedings{Liang2017,
    title={Why deep neural networks for function approximation?},
    author={S. Liang and R. Srikant},
    booktitle={International Conference on Learning Representations},
    year={2017}
}

@inproceedings{Lu2017,
    title = {The expressive power of neural networks: {A} view from the width},
    author = {Z. Lu and H. Pu and F. Wang and Z. Hu and L. Wang},
    booktitle = {Advances in Neural Information Processing Systems 30},
    year = {2017}

}

@inproceedings{Kingma2015,
  title={Adam: {A} method for stochastic optimization},
  author= {D. P. Kingma and J. Ba},
  booktitle={International Conference on Learning Representations},
  year={2015}
}

@article{Chui2018,
  title={Construction of neural networks for realization of localized deep learning},
  author={C. K. Chui and S. B. Lin and D. X. Zhou},
  journal={Front. Appl. Math. Stat.},
  volume={4},
  pages={14},
  year={2018},
}

@article{Guliyev2018,
  title={Approximation capability of two hidden layer feedforward neural networks with fixed weights},
  author={N. J. Guliyev and V. E. Ismailov},
  journal={Neurocomputing},
  volume={316},
  pages={262--269},
  year={2018}
}

@article{Maiorov1999,
  title={Lower bounds for approximation by {MLP} neural networks},
  author={V. Maiorov and A. Pinkus},
  journal={Neurocomputing},
  volume={25},
  pages={81--91},
  year={1999},
}

@article{Montanelli2019,
    title = {New error bounds for deep {ReLU} networks using sparse grids}, 
    author = {H. Montanelli and Q. Du},
    journal = {SIAM J. Math. Data Sci.},
    volume = {1},
    pages = {78--92},
    year = {2019}
}

@inproceedings{Suzuki2019,
    title={Adaptivity of deep {ReLU} network for learning in {B}esov and mixed smooth {B}esov spaces: optimal rate and curse of dimensionality},
    author={T. Suzuki},
    booktitle={International Conference on Learning Representations},
    year={2019}
}

@article{Yarotsky2017,
    title= {Error bounds for approximations with deep {ReLU} networks} ,
    author={D. Yarotsky},
    journal = {Neural Netw.},
    volume = {94},
    pages = {103--114},
    year = {2017}
}

@article{Petersen2018,
    title = {Optimal approximation of piecewise smooth functions using deep {ReLU} neural networks},
    author = {P. Petersen and F. Voigtlaender},
    journal = {Neural Netw.},
    volume = {108},
    pages = {296--330},
    year = {2018}
}

@article{Montanelli2021,
    title = {Deep {ReLU} networks overcome the curse of dimensionality for generalized bandlimited functions},
    author = {H. Montanelli and H. Yang and Q. Du},
    journal = {J. Comput. Math.},
    volume = {39},
    pages = {801--815},
    year = {2021}
}

@article{Shen2020,
   title = {Deep network approximation characterized by number of neurons},
   author = {Z. Shen and H. Yang and S. Zhang},
   journal = {Commun. Comput. Phys.},
   volume  = {28},
   pages = {1768--1811},
   year = {2020}
}

@article{Shen2019,
    title = {Nonlinear approximation via compositions},
    author = {Z. Shen and H. Yang and S. Zhang},
    journal = {Neural Netw.},
    volume = {119},
    pages = {74--84},
    year = {2019}
}

@inproceedings{Yarotsky2018,
    title={Optimal approximation of continuous functions by very deep {ReLU} networks},
    author={D. Yarotsky},
    booktitle={Proceedings of the 31st Conference on Learning Theory},
    year={2018}
}

@article{Chkifa2015,
  title={Breaking the curse of dimensionality in sparse polynomial approximation of parametric {PDE}s},
  author={A. Chkifa and A. Cohen and C. Schwab},
  journal={J. Math. Pures Appl.},
  volume={103},
  pages={400--428},
  year={2015}
}

@inproceedings{Parhi2023,
  title={Modulation spaces and the curse of dimensionality},
  author={R. Parhi and M. Unser},
  booktitle={2023 International Conference on Sampling Theory and Applications (SampTA)},
  year={2023}
}

@article{Zhi2020,
  title={Frequency principle: {F}ourier analysis sheds light on deep neural networks},
  author={Z. J. Xu and Y. Zhang and T. Luo and Y. Xiao and Z. Ma},
  journal={Commun. Comput. Phys.},
  volume={28},
  pages={1746--1767},
  year={2020}
}

@inproceedings{Xu2019,
  title={Training behavior of deep neural network in frequency domain},
  author={Z. J. Xu and Y. Zhang and Y. Xiao},
  booktitle={Neural Information Processing: 26th International Conference},
  year={2019}
}

@inproceedings{Cao2021,
  title={Towards understanding the spectral bias of deep learning},
  author={Y. Cao and Z. Fang and Y. Wu and D.-X. Zhou and Q. Gu},
  booktitle={Proceedings of the Thirtieth International Joint Conference on Artificial Intelligence},
  year={2021}
}

@inproceedings{Rahaman2019,
  title={On the spectral bias of neural networks},
  author={N. Rahaman and A. Baratin and D. Arpit and F. Draxler and M. Lin and F. Hamprecht and Y. Bengio and A. Courville},
  booktitle= {Proceedings of the 36th International Conference on Machine Learning},
  year={2019}
}

@article{Xu2024,
  title={Overview frequency principle/spectral bias in deep learning},
  author={Z. J. Xu and Y. Zhang and T. Luo},
  journal={Commun. Appl. Math. Comput.},
  volume={7},
  pages={827--864},
  year={2025}
}

@article{Luo2019,
  title={Theory of the frequency principle for general deep neural networks},
  author={T. Luo and Z. Ma and Z. J. Xu and Y. Zhang},
  journal={CSIAM Trans. Appl. Math.},
  volume={2},
  pages={484--507},
  year={2019}
}

@article{Luo2022,
  title={On the exact computation of linear frequency principle dynamics and its generalization},
  author={T. Luo and Z. Ma and Z. J. Xu and Y. Zhang},
  journal={SIAM J. Math. Data Sci.},
  volume={4},
  pages={1272--1292},
  year={2022}
}

@article{Zhang2021,
  title={A linear frequency principle model to understand the absence of overfitting in neural networks},
  author={Y. Zhang and T. Luo and Z. Ma and Z. J. Xu},
  journal={Chin. Phys. Lett.},
  volume={38},
  pages={038701},
  year={2021}
}

@inproceedings{Bordelon2020,
  title={Spectrum dependent learning curves in kernel regression and wide neural networks},
  author={B. Bordelon and A. Canatar and C. Pehlevan},
  booktitle={Proceedings of the 37th International Conference on Machine Learning},
  year={2020}
}

@inproceedings{Luo2022_2,
  title={An upper limit of decaying rate with respect to frequency in linear frequency principle model},
  author={T. Luo and Z. Ma and Z. Wang and Z. J. Xu and Y. Zhang},
  booktitle={Proceedings of Mathematical and Scientific Machine Learning},
  year={2022}
}

@article{Klusowski2018,
  title     = {Approximation by combinations of {ReLU} and squared {ReLU} ridge functions with $\ell^1$ and $\ell^0$ controls},
  author    = {J. M. Klusowski and A. R. Barron},
  journal = {IEEE Trans. Inform. Theory},
  volume    = {64},
  pages    = {7649--7656},
  year      = {2018}
}

@inproceedings{Bu2021,
    author = {J. Bu and Anuj A. Karpatne},
    title = {Quadratic Residual Networks: A New Class of Neural Networks for Solving Forward and Inverse Problems in Physics Involving PDEs},
    booktitle = {Proceedings of the 2021 SIAM International Conference on Data Mining},
    year={2021},
}

@inproceedings{Mantini2020,
  author={P. Mantini and S. K. Shah},
  title={CQNN: Convolutional Quadratic Neural Networks}, 
  booktitle={25th International Conference on Pattern Recognition},
  year={2020},
}

@article{Jiang2020,
  title={Nonlinear CNN: improving CNNs with quadratic convolutions},
  author={Y. Jiang and F. Yang and H. Zhu and D. Zhou and X. Zeng},
  journal={Neural Comput. Appl.},
  volume={32},
  pages={8507–-8516},
  year={2020}
}
\end{document}